\documentclass[paper]{siamart1116}

\usepackage{lipsum}
\usepackage{amsfonts}
\usepackage{graphicx}
\usepackage{epstopdf}
\def\NZQ{\Bbb}               

\def\RR{{\NZQ R}}
\usepackage{algorithmic}
\ifpdf
  \DeclareGraphicsExtensions{.eps,.pdf,.png,.jpg}
\else
  \DeclareGraphicsExtensions{.eps}
\fi

\numberwithin{theorem}{section}

\newcommand{\TheTitle}{Hybrid thermostatic approximations of junctions for some optimal control problems on networks} 
\newcommand{\TheAuthors}{F. Bagagiolo, and R. Maggistro}

\headers{hybrid approximations of junctions on networks}{\TheAuthors}

\title{{\TheTitle}\thanks{Submitted to the editors DATE.
\funding{This work was partially supported by OptHySYS project of the University of Trento.}}}

\author{
  Fabio Bagagiolo \thanks{Department of Mathematics, University of Trento, Via Sommarive 14, I-38123 Povo-Trento, Italy (\email{fabio.bagagiolo@unitn.it})}
  \and
  Rosario Maggistro \thanks{Department of Mathematical Sciences, Politecnico di Torino, Corso Duca degli Abruzzi, 24, I-10129 Torino, Italy (\email{rosario.maggistro@polito.it})}
}
\usepackage{amsopn}

\ifpdf
\hypersetup{
  pdftitle={\TheTitle},
  pdfauthor={\TheAuthors}
}
\fi

\begin{document}

\maketitle

\begin{abstract}
We study some optimal control problems on networks with junctions, approximate the junctions by a switching rule of delay-relay type and study the passage to the limit when $\varepsilon$, the parameter of the approximation, goes to zero. First, for a twofold junction problem we characterize the limit value function as viscosity solution and maximal subsolution of a suitable Hamilton-Jacobi problem. Then, for a threefold junction problem we consider two different approximations, recovering in both cases some uniqueness results in the sense of maximal subsolution. 
\end{abstract}



\section{Introduction}
\begin{figure}[htbp]
  \centering
\includegraphics[scale=0.42]{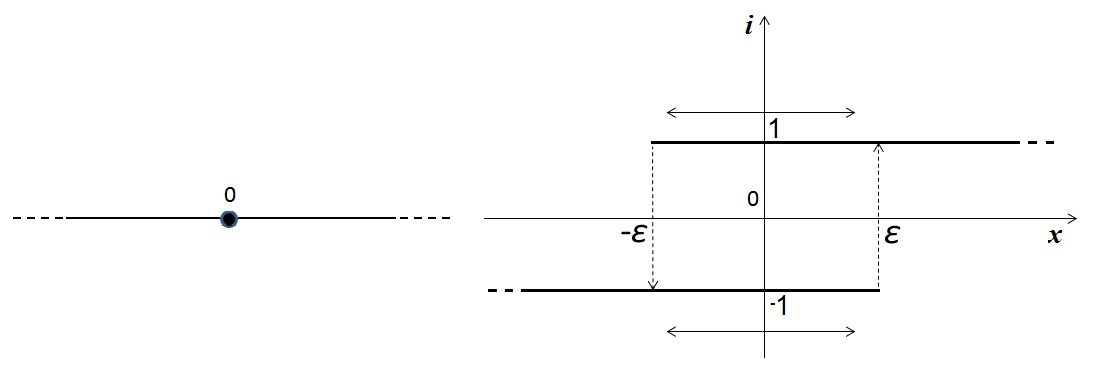}
  \caption{\label{fig:therm} The two-fold junction and its thermostatic approximation.}
  \end{figure}

In this paper, we are interested in optimal control problems with dynamics inside a network. Each arc $E_i$ of the network has its own controlled dynamics $f_i$ and cost $\ell_i$. When passing from an arc to another one through a node, the system then drastically experiences  a discontinuity.  We refer to this situation as a ``junction''. Recently there 
was an increasing interest in dynamical systems and differential equation on network, for example in connection with problems of data transmission and traffic flows (e. g. Garavello-Piccoli \cite{GaPi}, Engel et al. \cite{EnKrNa}).\\
For dynamic programming and Hamilton-Jacobi-Bellman (HJB) equation, in optimal control, the presence of junctions is a problem because, by the discontinuous feature of HJB, the uniqueness of the solution of HJB is not in general guaranteed. In particular for an optimal control problem we cannot characterize the value function as the unique solution of HJB. Some authors have recently studied optimal control and HJB on networks as well as discontinuous HJB also not necessarily coming from an optimal control problem, see for instance  Achdou-Camilli-Cutr\`\i-Tchou \cite{AcCaCuTc}, Camilli-Marchi \cite{CaMa}, Camilli-Marchi-Schieborn \cite{CaMaSc}, Camilli-Schieborn \cite{CaSc},  Imbert-Monneau-Zidani \cite{ImMoZi}, Achdou-Oudet-Tchou \cite{AcOuTc}, Achdou-Tchou \cite{AcTc}, and the recent Lions-Souganidis \cite{LiSo}. The optimal control problem on networks is related to $n$-dimensional optimal control problems on multi-domains, where the dynamics and costs incur in discontinuities when crossing some fixed hypersurfaces. These problems, started with Bressan-Hong \cite{BrHo}, \cite{BrHo1}, have been studied, in connection with HJB, in  Barles-Briani-Chasseigne \cite{BaBrCh}, Barnard-Wolenski \cite{BaWo}, Rao-Zidani \cite{RaZi},
Barles-Briani-Chasseigne \cite{BaBrCh2},  Rao-Siconolfi-Zidani \cite{RaSiZi}, Barles-Chasseigne \cite{BaCh}, Achdou-Oudet-Tchou \cite{achoudtch}, Barles-Briani-Chasseigne-Imbert \cite{BaBrChIm}, Imbert-Monneau\cite{ImMo}.\\
In this paper we study a possible approximation of an optimal control problem on a network with a junction. Some preliminary and partial results have been presented in Bagagiolo \cite{Bae}. The main critical point is a uniqueness result for the viscosity solution of HJB equation that turns out to be discontinuous in the (one-dimensional) space-variable (we refer the reader to, for example, Bardi-Capuzzo Dolcetta \cite{BaCaDo} for a comprehensive account of viscosity solutions for Hamilton-Jacobi equations). Indeed, when using the classical double-variable technique for proving comparison results between sub- and super- viscosity solutions we cannot in general conclude as in the standard way because the points of minimum and of maximum, even if very close, may belong to different arcs for which dynamics and costs are absolutely non-comparable (the junction, indeed). Consider the situation where two half-lines (the edges) are separated by one point (the junction). Also note that the possible angle between the lines is irrelevant, being the discontinuity of dynamics and costs through the junction the only relevant fact. Our approach is to replace the junction, which represents a unique threshold for passing both from one edge to the other one and vice-versa, by a so-called delayed thermostat consisting in two different thresholds for passing separately from one edge to the other one and vice-versa (see \cref{fig:therm}). The problem is then transformed in a so-called hybrid problem (continuous/discrete evolution, see for example Goebel-Sanfelice-Teel \cite{Goebel}) for which the discontinuity of HJB is replaced by some suitable mutually exchanged boundary conditions on the extreme points of the two branches. This allows to get a uniqueness result for HJB  for this kind of thermostatic problem. It is not unusual in the engineering literature on control problems to overcome discontinuities in dynamics by inserting some regularizing effects as hysteresis: switching and/or continuous. The thermostat is the fundamental brick of switching hysteresis model. See for instance Kokotovich \cite{TaoKoko} (pp. 17-23), Seidman \cite{Seidman88}, Hante-Leugering-Seidman \cite{HanteSeidman}. In Ceragioli-De Persis-Frasca \cite{Ceragioli} and Seidman \cite{Seidman2013} the thermostat is applied to solve several control problems coming from different contexts and discontinuities. The approximation of sliding mode behaviour through a thermostatic switching rule and the consequent passage to the limit of the switching threshold is discussed in Liberzon \cite{Liberzon} (pp. 14-15),  Utkin \cite{Utkin} (pp. 30-31) and Alexander-Seidman \cite{AlexanderSeidman}.
The study of this kind of switching hysteresis in connection with Dynamic Programming and HJB for optimal control is not so advanced and moreover the limit problem when the switching threshold goes to zero was probably never studied. Hence our results may shed light for new formulations for the junction problem as indeed it happens for our three-fold problem that we will study below. We start from the results in Bagagiolo \cite{Ba} (see also Bagagiolo-Danieli \cite{BaDa}), where the author studies the dynamic programming method and the corresponding HJB problem for optimal control problems whose dynamics has a thermostatic behavior. This means that the dynamics $f$ (as well as the cost $\ell$) besides the control, depends on the state variable $x \in \RR$, which evolves with continuity via the equation $x'=f$, and also depends on a discrete variable $i\in \left\{-1, 1\right\}$, whose evolution is governed by a delayed thermostatic rule, subject to the evolution of the state $x$. 
In \cref{fig:therm} the behavior of such a rule  is explained, correspondingly to the choice of the fixed threshold parameter $\varepsilon>0$. The output $i\in\{-1,1\}$ can (and must) jump  from $1$ to $-1$ only when the input $x$, coming from the right (i.e. from values larger than or equal to $-\varepsilon$), possibly goes below the threshold $-\varepsilon$; it can (and must) jump from $-1$ to $1$ only when $x$, coming from the left (i.e. from values smaller than or equal to $\varepsilon$) possibly goes above the threshold $\varepsilon$. In all other situations it remains constant. In particular, when $x>\varepsilon$ then $i$ is equal to $1$, and when $x<-\varepsilon$ then $i$ is equal to $-1$. We refer to Visintin \cite{Vi} page 102 for a formal definition of such a switching rule.
The controlled evolution is then given by
\[
\begin{cases}
x'(t) = f(x(t), i(t), \alpha(t)), \\
i(t)= h_{\varepsilon}\left[x\right](t) \\
x(0)= x_0, \ \ \ i(0)= i_0
\end{cases} \]
where $\alpha: \left[0, +\infty\right] \rightarrow A$ is the measurable control, and $h_{\varepsilon}\left[\cdot\right]$ represents the thermostatic delayed relationship between the input $x$ and the output $i$. The initial value $i_0 \in \left\{-1, 1\right\}$ must be coherent with the thermostatic relation: $i_0 = 1$ (resp. $i_0 = -1$) whenever $x_0 > \varepsilon$ (resp. $x_0 < - \varepsilon$). The infinite horizon optimal control problem is then, given a running cost $\ell$ and a discount factor $\lambda > 0$, the minimization over all measurable controls, of the cost
\begin{equation}\label{eq:[4]}
V_{\varepsilon}(x_0, i_0)= \inf_{\alpha \in{\cal A}} \int_{0}^{\infty} e^{-\lambda t} \ell(x(t), i(t), \alpha(t))dt 
\end{equation}
where ${\cal A}$ is the set of measurable controls.
In \cite{Ba}, the problem is written as a coupling of two exit time optimal control problems which mutually exchange their exit-costs. In particular, using the notations
\begin{equation}
\label{eq:omega}
\Omega_1^{\varepsilon}= \left\{x > -\varepsilon \right\},\ \ \overline\Omega_1^{\varepsilon}= \left\{x \geq -\varepsilon \right\},\ \ \Omega_{-1}^{\varepsilon}= \left\{x < \varepsilon \right\} \ \ \overline\Omega_{-1}^{\varepsilon}= \left\{x \leq \varepsilon \right\},
\end{equation}
in $\overline\Omega_1^{\varepsilon} \times \left\{1\right\}$ (resp. $\overline\Omega_{-1}^{\varepsilon} \times \left\{-1\right\}$), the function $x\mapsto V_{\varepsilon}(x,1)$ (resp. $x\mapsto V_\varepsilon(x,-1)$) coincides with the value function of the exit-time optimal control problem on $\overline\Omega_1^{\varepsilon}$ \ (resp. $\overline\Omega_{-1}^{\varepsilon}$ ), where the exit-cost on $-\varepsilon$\ (resp. on $\varepsilon$) is given by $V_\varepsilon(-\varepsilon,-1)$ \ (resp. $V_\varepsilon(\varepsilon,1)$). Under standard hypotheses, in \cite{Ba} is proved the following theorem.
\begin{theorem}\label{th:sistemapprox}
The value function $V_{\varepsilon}$ in \cref{eq:[4]} is the unique bounded, continuous viscosity solution on $\overline\Omega_1^{\varepsilon} \times \left\{1\right\} \cup \overline\Omega_{-1}^{\varepsilon} \times \left\{-1\right\}$ of the following coupled Dirichlet problem, where the boundary conditions (the two exit-costs) are also in the viscosity sense ($V'$ stays for derivative with respect to $x$)
\begin{equation}\label{s1}
\begin{cases}
\lambda V_{\varepsilon}(x, 1) + \sup_{a \in A}\left\{-f(x, 1, a)V'_{\varepsilon}(x, 1) - \ell(x, 1,a)\right\} = 0 \ \text{in}\ \Omega_1^{\varepsilon} \times \left\{1\right\} \\
V_{\varepsilon}(- \varepsilon, 1)= V_{\varepsilon}(- \varepsilon, -1) \\
\lambda V_{\varepsilon}(x, -1) + \sup_{a \in A}\left\{-f(x, -1, a)V'_{\varepsilon}(x, -1) - \ell(x, -1, a)\right\} = 0 \ \text{in}\ \Omega_{-1}^{\varepsilon} \times \left\{-1\right\} \\
V_{\varepsilon}(\varepsilon, -1)= V_{\varepsilon}(\varepsilon, 1)
\end{cases} 
\end{equation}
\end{theorem}

In the present paper, we approximate some junction problems by a suitable combinations of delayed thermostats. Every thermostat is characterized by its two thresholds ($\varepsilon$ and $-\varepsilon$ in the preceding description). We study the limit of the value functions $V_\varepsilon$ and of the HJB problem when the threshold distance $\varepsilon$ tends to zero, and hence recovering the junction situation.

\begin{figure}[htbp]
\begin{center}
\includegraphics[width=%
0.45\textwidth]{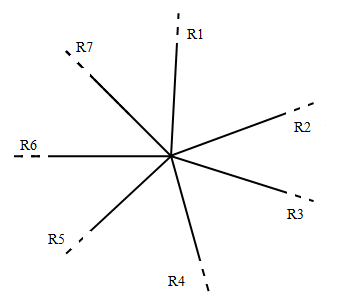} 
\caption{\label{7foldjunction} A star-shaped network (7-fold junction).}
\end{center}
\end{figure}

In Barles-Briani-Chasseigne \cite{BaBrCh}, among others, a one-dimensional two-fold junction problem is studied and some possible approximations are given. Here we introduce a different kind of approximation (thermostatic) and recover, by a different proof, similar results: we characterize the limit problem and we get that the limit of $V_\varepsilon$ is the corresponding maximal viscosity subsolution. This corresponds to the value function of the junction optimal control problem where, on the junction point, some further dynamics (besides the already given ones) are considered: the ones given by a suitable convexification of ``non-inward pointing" dynamics (``regular" dynamics in \cite{BaBrCh}) and somehow corresponding to stable equilibria on the junction point (stable equilibria of dynamics interpreted as forces).
In \cite{BaBrCh} the case where on the junction one can also use the so-called ``singular" dynamics (i.e., a suitable convexification of ``inward pointing" dynamics; somehow corresponding to unstable equilibria) 
is also treated. In particular in this case all possible dynamics on the junction can be used: singular, regular and the already given ones of the original control problem.
All such possible dynamics at the interface are also used in 
Rao-Siconolfi-Zidani \cite{RaSiZi} to prove uniqueness of the solution for HJB.
In our paper, the use of the thermostatic approximation leads to consider only the ``regular problem", and hence to have a characterization as maximal subsolution. Hence the problem in \cite{RaSiZi} and our limit problem are substantially different.
Indeed, as said before, the concept of thermostat is
based on the concept of switching that occurs when suitable thresholds are reached. For the occurrence of the switching it is necessary that the threshold is reached with a suitable signed velocity.
And such a sign is exactly the one requested by the construction of the regular dynamics.
 Moreover, it is important to note that there exist several ways to define a junction optimal control problem and everyone of them has different HJ representation with different possible approximations by more regular problem.
The case of a ``threefold" junction (see \cref{threejunction}) is not treated in \cite{BaBrCh}, and indeed the convexification of dynamics seems to be not more applicable (the physical interpretation as forces equilibrium is also failing). However, inspired by the previous thermostatic approximation, we introduce a special kind of ``convexification parameters" that somehow corresponds to the length of the time intervals that the dynamics spends on every single branch of a ``threefold" thermostatic approximation.  Here, we have more than one way for passing to the limit, and we recover uniqueness results for the limit problems in the sense of maximal subsolution. The approach can be extended to the $n$-fold junctions (\cref{7foldjunction}), but with higher notation complexity.

The paper is organized as follows: basic assumption are set in \cref{sect:2}. In \cref{sec:3}, we introduce the thermostatic approximation of a two-junction problem and study the passage to the limit for $\varepsilon \rightarrow 0$ in the problem studied in \cite{Ba}. In \cref{sec:4} we study the three-junction problem both in the case with uniform switching thresholds and with non uniform switching thresholds. This corresponds to two different limit optimal control problems with different admissible behaviors on the junction.

\section{Basic assumption on the junction problem and the delay thermostat}
\label{sect:2}
Let the junction be given by a finite number of co-planar half-lines $R_i$, $i=1,\dots,n$, originating from the same point $O$, and we consider the half-lines as closed, that is the point $O\in R_i$ for every $i$. On each branch $R_i$ we consider a one-dimensional coordinate $x\ge0$ such that $x(O)=0$. The state position may be then encoded by the pair $(x,i)$. In Sect. 4, for convenience of notation, we will sometimes change the sign of $x$ .
We consider a controlled evolution on such a star-shaped network, given by the following dynamics. 
On $R_i$ the system is driven by a continuous and bounded dynamics $f_i:\RR\times A\to\RR$, with $A$ compact, $f_i$ Lipschitz and controllability holds:

\begin{equation}\label{eq:Lip}
\exists L>0\ \mbox{such that } \forall\ x, y \in \RR,\ \forall\ a \in A\ \mbox{it is }
\vert f_i(x, a)- f_i(y, a) \vert \leq L \vert x-y \vert.
\end{equation}

\begin{equation}\label{eq:Controllability}
\forall\ i\ \exists \ a_{i}^{-}, a_{i}^{+} \ \in A\ \ \text{s.t.}\ \ \ f_{i}(0, a_{i}^{-}) < 0 < f_{i}(0, a_{i}^{+})  
\end{equation}
The controlled system on the network is then, for an initial state $(x,i)$ with $x\in R_i$,
\begin{equation}
\label{eq:systemjunction}
\left\{
\begin{array}{ll}
\displaystyle
y'(t)=f_j(y(t),\alpha(t))&\mbox{for } t>0\ \mbox{and } y(t)\in R_j\\
y(0)=x\\
x\in R_i
\end{array}
\right.
\end{equation}

\noindent
where $\alpha:[0,+\infty[\to A$ belongs to $\cal A$, the set of measurable controls, and $j=j(t)$ is the switching variable that switches to $j'$ when $y(t)$ enters the new half-line $R_{j'}$.

To this controlled systems we associate an infinite horizon optimal control problem. For every branch $R_i$ we consider a running cost $\ell_i:\RR\times A\to [0,+\infty[$, and the problem is given by the minimization, over all controls $\alpha\in{\cal A}$, of the cost functional
\begin{equation}
\label{eq:costfunctionaljunction}
J(x,i,\alpha)=\int_0^{+\infty}e^{-\lambda t}\ell_j(y(t),\alpha(t))dt.
\end{equation}
\noindent
In \cref{eq:costfunctionaljunction}, $\lambda>0$ is a fixed discount factor, the trajectory $y(\cdot)$ is the solution of \cref{eq:systemjunction}, and the index $j$ switches as explained above. Moreover, for every $i$, the function $\ell_i:\RR \times A\rightarrow \RR$ is continuous and bounded, and there exists a modulus of continuity $\omega_{\ell}:[0,\infty[\to[0,+\infty[$ (i.e. continuous, increasing and $\omega_\ell(0)=0$), such that, for any $x, y \in \RR$ and $a \in A$ and for any $i$
\begin{equation}\label{eq:LLip}
\vert \ell_i(x, a)- \ell_i(y, a) \vert \leq \omega_{\ell}\left( \vert x-y \vert\right).
\end{equation}

We finally consider the value function
\begin{equation*}
V(x,i)=\inf_{\alpha\in{\cal A}}J(x,i,\alpha).
\end{equation*}

Of course, the concept of solution (or trajectory) for the system \cref{eq:systemjunction} and the definition of the cost \cref{eq:costfunctionaljunction} are not well-posed. Indeed, at the junction point $O$, we can choose the index $i$ we prefer, but the existence of the trajectory is not guaranteed, due to possible fast oscillations of the index $i$ (as when, for a generic ordinary differential equation, the dynamics is discontinuous in the space variable). The main goal of the present paper is just an approximation and the corresponding passage to the limit, for such possible oscillating behavior in the context of optimal control. To this end, we are going to use delayed thermostat operator and, in addition to what already said in the Introduction, here we point out that, fixed the thresholds $-\varepsilon,\varepsilon$, for each continuous scalar input $t\mapsto u(t)$, and for each initial output $i_0\in\{-1,1\}$ coherent with $u(0)$, there exists a unique output $t\mapsto i(t)=:h_\varepsilon[u](t)\in\{-1,1\}$ satisfying $i(0)=i_0$. For a regular scalar dynamics $g$, and for a coherent initial state $(x,i_0)$, there exists a unique solution of the thermostatic system

\[
\left\{
\begin{array}{ll}
y'=g(y,i,t)\\
y(0)=x\\
i(t)=h_\varepsilon[y](t),\ \ i(0)=i_0
\end{array}
\right.
\]

\noindent
The main  reason for that is indeed the ``splitting" of the thresholds, which avoids fast oscillations of the switching variable $i$, and allows to construct the solution by a suitable gluing of pieces of solutions with constant $i$ (see \cite{Ba}).

\section{A twofold junction problem}\label{sec:3}
We consider a one-dimensional optimal control problem for which the controlled dynamics and the cost, $f,\ell$, suddenly change when passing from one half-line to the other one: $f(x, \cdot) = f_1(x, \cdot)\ (\text{resp.}\ f(x, \cdot) = f_{-1}(x, \cdot))$ if $x > 0 \ (\text{resp. if}\ x < 0)$, where $f_1 : [0, + \infty[ \times A \rightarrow \RR$, $f_{-1} : ] -\infty, 0] \times A \rightarrow \RR$. The point $x=0$ may represent a ``junction'', a node on a network with two entering edges (see \cref{fig:therm}). For $\varepsilon > 0$ we approximate the junction problem by a delayed thermostatic problem. Still denoting by $f_1, f_{-1}$ two extensions by constancy in the space variable $x$ of the dynamics to $[-\varepsilon, +\infty[ \times A$ and to $]-\infty, \varepsilon] \times A$ respectively, we may consider the controlled system
\begin{equation}\label{s2}
\begin{cases}
x'(t) = f_{i(t)}(x(t), \alpha(t)), \\
i(t)= h_{\varepsilon}\left[x\right](t) \\
x(0)= x_0, \ \ \ i(0)= i_0
\end{cases} 
\end{equation}
Similarly to $f_1,f_{-1}$, we extend the running costs $\ell_1, \ell_{-1}$. 

Let $V_{\varepsilon}$ be the value function of the thermostatic optimal control problem with dynamics given by \cref{s2} and corresponding costs. We define the function
\begin{flalign*}
\tilde{V}_{\varepsilon} : \RR \setminus \left\{0\right\} \rightarrow \RR,  & \ \
\tilde{V}_{\varepsilon} (x)= 
\begin{cases}
V_{\varepsilon}(x, 1) & x > 0\\
V_{\varepsilon}(x, -1) & x < 0.
\end{cases}
\end{flalign*}
In general, $V_{\varepsilon}(0, -1) \neq V_{\varepsilon}(0, 1)$.
\begin{theorem} \label{th:2junction}
As $\varepsilon \rightarrow 0^{+}$, $\tilde{V}_{\varepsilon}$ uniformly converges on $\RR\setminus  \left\{0\right\}$ to a continuous function which, if \cref{eq:Lip,eq:LLip} hold, continuously extends to a function $\tilde V$ on the whole $\RR$. If \cref{eq:Controllability} holds, $\tilde V$ coincides with the (already known as unique) viscosity solution of 
\begin{equation}\label{eq:Sistema2junction}
\begin{cases}
\lambda V + H_{1}(x, V') = 0 \ \ \ \text{for} \  x>0\\
\lambda V + H_{-1}(x, V') = 0 \ \ \text{for} \  x<0 \\
V(0) = \min \left\{u_{0}(0), V_{sc(-1)}(0), V_{sc(1)}(0) \right\}
\end{cases}
\end{equation}
where $H_1,H_{-1}$ are the Hamiltonians in \cref{s1}, $u_0(0)$ is the convexification
\begin{equation} \label{eq:date}
u_0(0)= \frac{1}{\lambda}\min_{A_0}\left\{\mu \ell_{-1}(0, a_{-1}) + (1 - \mu)\ell_{1}(0, a_{1})\right\}
\end{equation}
\begin{multline}
\label{eq:A0}
A_0 = \bigl\{(\mu, a_{-1}, a_{1} ) \in [0,1] \times A \times A :\\
 \mu f_{-1}(0, a_{-1}) + (1 - \mu)f_{1}(0, a_{1})=0, f_1(0, a_1) \leq 0, f_{-1}(0, a_{-1})\geq 0 \bigr\}
\end{multline}
and $V_{sc(i)}(0)$ is the value function at $x=0$ of the state-constraint optimal control problem on the branch $i$.
\end{theorem}
\begin{proof}
We are going to use the notations in \cref{eq:omega}. We also recall that the state-constraint problem in branch $i$ is the optimal control problem restricted to the branch $i$ and such that, at the point $x=0$ we can only use controls that make us to not leave the  branch. We first prove that $V_\varepsilon$ uniformly converges to a continuous function on $\mathbb{R}\setminus\{0\}$. We have some cases and we illustrate some of them. 

i) $f_{-1}(0,a)\le0$ for all $a\in A$. Hence, when starting from a point of $\overline\Omega_{-1}^{\varepsilon} \times \left\{-1\right\}$, it is impossible to switch on the other branch $\overline\Omega_1^{\varepsilon} \times \left\{1\right\}$. Hence, $x\to V_\varepsilon(x,-1)$ is the value function of the optimal control problem with dynamics $f_{-1}$ and cost $\ell_{-1}$ and state-constraint in $\overline\Omega_{-1}^{\varepsilon} \times \left\{-1\right\}$, which uniformly converges on $]-\infty,0]$ to the value function with same dynamics and cost and with state constraints in $]-\infty,0]$ (dynamics and costs are bounded), that is to $V_{sc(-1)}$. In the other branch, being $V_{\varepsilon}(-\varepsilon,-1)$ convergent to $V_{sc(-1)}(0)$, we also get the uniform convergence of $V_\varepsilon(\cdot,1)$ to the unique solution of the first line of \cref{eq:Sistema2junction} with viscosity boundary datum $V_{sc(-1)}(0)$. Indeed, they are respectively the value function of the exit-time problem in $[-\varepsilon,+\infty[$ and $[0,+\infty[$ with the same dynamics, same cost and with convergent exit-costs.

ii) $\exists \ a_{-1},a_1\in A$ such that $f_1(0,a_1)<0<f_{-1}(0,a_{-1})$. In this case, when $\varepsilon$ is sufficiently small, starting from $(\varepsilon,1)$ (resp. from $(-\varepsilon,-1)$) we can always switch on the other branch, and we can reach  $(-\varepsilon,-1)$ (resp. $(\varepsilon,1)$) in a time interval whose length is infinitesimal as $\varepsilon$. It is then easy to check that the difference $|V_\varepsilon(\varepsilon,1)-V_\varepsilon(-\varepsilon,-1)|$ (as well as $|V_\varepsilon(0,1)-V_\varepsilon(0,-1)|$) is also infinitesimal as $\varepsilon$. 
Moreover, for every pair $(\varepsilon_1,\varepsilon_2)$ with $\varepsilon_1,\varepsilon_2>0$, $\|V_{\varepsilon_1}-V_{\varepsilon_2}\|$ is also infinitesimal as $\max\{\varepsilon_1,\varepsilon_2\}$. As before, $V_\varepsilon$ uniformly converges on $\mathbb{R}\setminus\{0\}$ to a solution of the first two lines of \cref{eq:Sistema2junction}, which also continuously extends to $x=0$.
We denote by $\tilde V$ such extended limit function and, assuming \cref{eq:Controllability} (which of course implies the conditions in ii)), we prove that it satisfies the third equation of \cref{eq:Sistema2junction}, from which the conclusion of the proof. Again, we proceed illustrating some cases.

a) $V_{sc(-1)}(0)$ strictly realizes the minimum in \cref{eq:Sistema2junction}. Then, there exists a measurable control $\alpha$ such that the corresponding trajectory starting from $x=0$ with dynamics $f_{-1}$ does not exit from $]-\infty,0]$, and the corresponding cost, with running cost $\ell_{-1}$, is strictly less than $u_{0}(0)$ and $V_{sc(1)}(0)$. The control $\alpha$ has exactly the same cost for the thermostatic problem with initial point $(0,-1)$ (no switchings occur).\\ 
Now, we observe that, for every $(\mu,a_{-1},a_1)\in A_0$ with $f_{-1}(0,a_{-1}),f_1(0,a_1)\neq0$, the alternation of the constant controls $a_{-1},a_1$ correspondingly to every switching, gives a cost for the thermostatic problem in $(0,-1)$ as well as in $(0,1)$, which, when $\varepsilon$ goes to zero, converges to $\left(\mu\ell_{-1}(0,a_{-1})+(1-\mu)\ell_1(0,a_1)\right)/\lambda$. Indeed, the condition $\mu f_{-1}(0,a_{-1})+(1-\mu)f_1(0,a_1)=0$ implies that $f_{-1}(0,a_{-1}),f_1(0,a_1)$ are in the same (inverse) proportion as $\mu$ and $1-\mu$, and the corresponding time-durations for covering the distance $2\varepsilon$ (from one threshold to the other one) are in the same (direct) proportion as $\mu$ and $1-\mu$. Hence, the required convergence holds.\\
From this we get that $V_{sc(-1)}(0)=V_\varepsilon(0,-1)$ and so $\tilde V(0)= V_{sc(-1)}(0)$.

b) $u_{0}(0)$ strictly realizes the minimum. Then, let $(\mu,a_{-1},a_1)\in A_0$ be such that $\mu\ell_{-1}(0,a_{-1})+(1-\mu)\ell_{1}(0,a_{1})$ is the minimum in the definition of $u_{0}(0)$. Again, as in the previous point, we get that a switching trajectory using controls $a_{-1}$ and $a_1$ is near optimal for $V_\varepsilon$, and then the conclusion.
\end{proof}
\newsiamremark{rem}{Remark}
\begin{rem}
In this one-dimensional case, \cref{th:2junction} also proves that $\tilde V=U^+$, where $U^+$ is the value function of the so-called regular problem in \cite{BaBrCh}. In the sequel we are also given a different proof of such an equality where, using the thermostatic approximation, we show that $\tilde V$ is the maximal subsolution of a suitable Hamilton-Jacobi problem as in \cite{BaBrCh}, namely next problem \cref{eq:HJBproblem}.
\end{rem}
\begin{theorem} \label{thm:ishii}
Assume \cref{eq:Lip,eq:Controllability,eq:LLip}. The function $\tilde{V}$ is a viscosity solution of the Hamilton-Jacobi-Bellman problem
\begin{equation}\label{eq:HJBproblem}
\begin{cases}
\lambda V + H_{1}(x, V') = 0 \ \ \text{in} \  \{x>0\}=:\Omega_1 \\
\lambda V + H_{-1}(x, V') = 0 \ \ \text{in} \  \{x<0\}=:\Omega_{-1} \\
\min \left\{\lambda V + H_1, \lambda V + H_{-1}\right\} \leq 0 \ \text{on} \ x=0 \\
\max \left\{\lambda V + H_1, \lambda V + H_{-1}\right\} \geq 0 \ \text{on} \ x=0.
\end{cases}
\end{equation}
\noindent
Here we mean that $\tilde V$ is a subsolution of the first three equations and a supersolution of the first two together with the fourth one.
\end{theorem}
\begin{proof}
From \cref{th:2junction}, $\tilde V$ is a viscosity solution of the first two lines of \cref{eq:HJBproblem}.

We now prove the third equation in  \cref{eq:HJBproblem}. Let $\varphi \in C^1(\RR)$ be a test function such that $\tilde V-\varphi$ has a strict relative maximum at $x=0$. 
By uniform convergence, there exists a sequence $x_\varepsilon\in\overline\Omega_1^\varepsilon$ of points of relative maxima for $V_\varepsilon(\cdot,1)-\varphi$ which converge to $x=0$. We may have two mutually exclusive cases: 1) at least for a subsequence,  at $x_\varepsilon$ the HJB equation satisfied by $V_\varepsilon(\cdot,1)$ has the right sign "$\le$" (if $x_\varepsilon$ is an interior point, then we always have the right sign, being the equation satisfied), 2) it is definitely true that the boundary point $x_\varepsilon=-\varepsilon$ is a strict maximum point and the HJB equation has the wrong sign "$>$". Also note that the boundary of $\overline\Omega_1^\varepsilon$, i.e. $x=-\varepsilon$, is also converging to $x=0$.

Case 1). As $\varepsilon\to0$, we get $\lambda \tilde{V}+H_1 \leq 0$ in $x=0$ and the third equation in \cref{eq:HJBproblem}. 

Case 2). Since the boundary conditions in \cref{s1} are in the viscosity sense and by virtue of the controllability condition \cref{eq:Controllability}, we have
\begin{equation} \label{eq:CondBordosup}
V_{\varepsilon}(- \varepsilon, 1) = V_{\varepsilon}(- \varepsilon, -1)
\end{equation}

The same argumentations and cases also hold for the branches $\overline\Omega_{-1}^\varepsilon$. If the corresponding case 1) holds, then we get the conclusion as before. Otherwise we have
\begin{equation} \label{eq:CondBordoinf}
V_{\varepsilon}(\varepsilon, -1) = V_{\varepsilon}(\varepsilon, 1)
\end{equation}

We prove that case 2) cannot simultaneously holds in both branches. Indeed, observe that $(- \varepsilon, -1)$ is in the interior of $\overline{\Omega}_{-1}^{\varepsilon}$ and $(\varepsilon,1)$ is in the interior of $\overline\Omega_1^\varepsilon$,  therefore, using 
\cref{eq:CondBordoinf,eq:CondBordosup}, we get the following contradiction and conclude
\begin{equation*}
\begin{array}{ll}
\displaystyle
V_{\varepsilon}(- \varepsilon, -1) - \varphi (- \varepsilon) < V_{\varepsilon}(\varepsilon, -1) - \varphi (\varepsilon )= V_{\varepsilon}(\varepsilon, 1) - \varphi (\varepsilon )\\
\displaystyle
< V_{\varepsilon}(- \varepsilon, 1) - \varphi (-\varepsilon) = V_{\varepsilon}(- \varepsilon, -1) - \varphi (- \varepsilon)
\end{array}
\end{equation*}
To prove the fourth equation in  \cref{eq:HJBproblem}, we argue in the same way.
\end{proof}
We prove that $\tilde{V}$ is the maximal subsolution of \cref{eq:HJBproblem}, and use the following lemma.

\begin{lemma}
\label{lem:vbarra}
Assume that $\forall\ \varepsilon >0$ small enough, the optimal strategy for the approximating problem $\varepsilon$, starting by any $(x, 1), (x, -1)$ with $x \in [-\varepsilon, \varepsilon]$, is  to have infinitely many switches between the two branches (i.e. no state-constraint behavior is optimal). Let
$(\bar{\mu}, \bar{a}_{-1}, \bar{a}_1) \in A_0$ be such that
$f_{-1}(0, \bar{a}_{-1})>0, f_1(0, \bar{a}_1)<0$, and that
\begin{equation}\label{eq:Optimality}
\tilde V(0)=u_0(0)= \frac{1}{\lambda}\lbrace \bar{\mu}\ell_{-1}(0, \bar{a}_{-1})+(1-\bar{\mu})\ell_1(0, \bar{a}_1) \rbrace.
\end{equation}
For every $x \in [-\varepsilon, \varepsilon]$, we consider the two switching trajectories (compare with \cref{s2})
\[
\begin{cases}
y'(t) = f_{i(t)}(0, \bar{a}_{i(t)}), \\
i(t)= h_{\varepsilon}\left[y\right](t), \\
y(0)= x, \ \ \ i(0)= 1,
\end{cases}, \ \ \  
\begin{cases}
y'(t) = f_{i(t)}(0, \bar{a}_{i(t)}), \\
i(t)= h_{\varepsilon}\left[y\right](t), \\
y(0)= x, \ \ \ i(0)= -1.
\end{cases}
\]
\noindent
On the branches, they have constant velocity ($f_1(0, \bar{a}_1)$ and $f_{-1}(0, \bar{a}_{-1})$ towards left and right respectively), and switch infinitely many times.
We consider the functions
\begin{equation}\label{vbarrainfinito}
\begin{aligned}
\bar{V}_{\varepsilon}(x, 1) & = \int_{0}^{\infty} e^{-\lambda t}\ell_{i(t)}(0, \bar{a}_{i(t)})dt \quad \text{with}\quad i(0)=1,\\
\bar{V}_{\varepsilon}(x, -1) & = \int_{0}^{\infty} e^{-\lambda t}\ell_{i(t)}(0, \bar{a}_{i(t)})dt \quad \text{with}\quad i(0)=-1.
\end{aligned}
\end{equation}
Then $\bar{V}_{\varepsilon}(\cdot, 1)$ and $\bar{V}_{\varepsilon}(\cdot, -1)$ are differentiable in $[-\varepsilon, \varepsilon]$ and
\begin{equation}\label{DifferenzaNulla}
\sup_{x \in [-\varepsilon, \varepsilon]}\vert \bar{V}_{\varepsilon}^{'}(x, 1)- \bar{V}_{\varepsilon}^{'}(x, -1) \vert \rightarrow 0\ \ \text{for} \ \ \varepsilon\rightarrow 0.
\end{equation}
\end{lemma}
\begin{proof} 
The derivability comes form the constancy of dynamics and costs. We can rewrite the two functions in \cref{vbarrainfinito} as 
\begin{equation}
\label{eq:vbarra}
\begin{aligned}
\bar{V}_{\varepsilon}(x, 1) & = \int_{0}^{\frac{x+\varepsilon}{\vert f_1(0, \bar{a}_1)\vert}} e^{-\lambda t}\ell_1(0, \bar{a}_1)dt + e^{\frac{-\lambda(x+\varepsilon)}{\vert f_1(0, \bar{a}_1)\vert}}\bar{V}_{\varepsilon}(-\varepsilon, 1),\\
\bar{V}_{\varepsilon}(x, -1) & = \int_{0}^{\frac{\varepsilon - x}{f_{-1}(0, \bar{a}_{-1})}} e^{-\lambda t}\ell_{-1}(0, \bar{a}_{-1})dt + e^{\frac{-\lambda(\varepsilon - x)}{ f_{-1}(0, \bar{a}_{-1})}}\bar{V}_{\varepsilon}(\varepsilon, -1),
\end{aligned}
\end{equation}
where the upper extremal of the integration is the reaching time of the threshold in the corresponding initial branch. Then we have
\begin{equation}\label{uguaglianzavbarra}
\bar{V}_\varepsilon(-\varepsilon, 1)= \bar{V}_\varepsilon(-\varepsilon, -1)\quad \text{and} \quad \bar{V}_\varepsilon(\varepsilon, -1)= \bar{V}_\varepsilon(\varepsilon, 1),
\end{equation}
and by \cref{eq:Optimality} for any $i$, $\lim_{\varepsilon\to 0}\bar{V}_\varepsilon(x,i)= \tilde{V}(0)= u_0(0)$ uniformly in $x \in [-\varepsilon, \varepsilon]$.
A direct calculation gives
\begin{align*}
\bar{V}_{\varepsilon}^{'}(x, 1)& = \frac{1}{\vert f_1(0, \bar{a}_1)\vert}e^{\frac{-\lambda(x+\varepsilon)}{\vert f_1(0, \bar{a}_1)\vert}}\ell_1(0, \bar{a}_1)-\frac{\lambda e^{\frac{-\lambda(x+\varepsilon)}{\vert f_1(0, \bar{a}_1)\vert}}}{\vert f_1(0, \bar{a}_1)\vert}\bar{V}_{\varepsilon}(-\varepsilon, 1),\\
\bar{V}_{\varepsilon}^{'}(x, -1)& = - \frac{1}{f_{-1}(0, \bar{a}_{-1})}e^{\frac{-\lambda(\varepsilon-x)}{ f_{-1}(0, \bar{a}_{-1})}}\ell_{-1}(0, \bar{a}_{-1})+\frac{\lambda e^{\frac{-\lambda(\varepsilon-x)}{ f_{-1}(0, \bar{a}_{-1})}}}{f_{-1}(0, \bar{a}_{-1})}\bar{V}_{\varepsilon}(\varepsilon, -1).
\end{align*}
and then for $\varepsilon \to 0$
\begin{align*}
& \bar{V}_{\varepsilon}^{'}(x, 1) \longrightarrow \frac{\bar{\mu}\big(\ell_1(0, \bar{a}_1)- \ell_{-1}(0, \bar{a}_{-1})\big)}{\vert f_1(0, \bar{a}_1)\vert}, \\
& \bar{V}_{\varepsilon}^{'}(x, -1) \longrightarrow \frac{(\bar{\mu}-1)\big(\ell_{-1}(0, \bar{a}_{-1})- \ell_{1}(0, \bar{a}_{1})\big)}{ f_{-1}(0, \bar{a}_{-1})}.
\end{align*}
\noindent
Recalling the definition of $A_0$ \cref{eq:A0}, calculating $\overline\mu$, we get \cref{DifferenzaNulla} by
\begin{equation*}
 \bar{V}_{\varepsilon}^{'}(x, 1) \longrightarrow \frac{\ell_{1}(0, \bar{a}_{1})- \ell_{-1}(0, \bar{a}_{-1})}{f_{-1}(0, \bar{a}_{-1})- f_1(0, \bar{a}_1) }, \quad
 \bar{V}_{\varepsilon}^{'}(x, -1) \longrightarrow 
\frac{\ell_{1}(0, \bar{a}_{1})- \ell_{-1}(0, \bar{a}_{-1})}{f_{-1}(0, \bar{a}_{-1})- f_1(0, \bar{a}_1)}
\end{equation*}
\end{proof}
\begin{theorem}\label{Confronto} For each $u$ bounded, continuous subsolution of \cref{eq:HJBproblem}, it is $u \leq \tilde{V}$ in $\RR$.
\end{theorem}
\begin{proof}
We can assume to be in the situation as in \cref{lem:vbarra}. Indeed, otherwise in at least one branch $\tilde V$ coincides with the corresponding state-constraint value function which, see for example Soner \cite{So}, is greater than any subsolution (note that, in general the state-constraint value functions do not satisfy the third line of \cref{eq:HJBproblem}). We then also get $u\leq \tilde{V}$ on the other branch. \\
We assume by contradiction that $\sup_{x \in \RR}(u-\tilde{V})(x)>\delta>0$. If 
\begin{equation*} 
\exists r>0 \vert \forall \delta' >0\ \exists\ \overline{x} \in ]r, +\infty[: \sup_{x\in \RR} \big((u-\tilde{V})(x)-(u-\tilde{V})(\overline{x})\big)\leq \delta',
\end{equation*}
then, by \cref{thm:ishii} and known comparison techniques we get a contradiction because, in $ ]r, +\infty[$, $\tilde V$ is a supersolution and $u$ is a subsolution of the same HJB. Similarly for the opposite case $]-\infty, -r[$. 
Hence we may restrict to the case where $u -\tilde{V}$ has the maximum with respect to $r$ in $x=0$.
Since $\bar{V}_{\varepsilon}(x, i)$ converges to $\tilde{V}(0)$, with $\bar V_\varepsilon$ defined in \cref{eq:vbarra}, then for small $\varepsilon$,
\begin{equation}\label{eq:ConditionMassimo}
u(z^i)-\bar{V}_{\varepsilon}(z^i, i) =\max_{[-\varepsilon, \varepsilon]}(u(\cdot)-\bar{V}_{\varepsilon}(\cdot, i))> \frac{\delta}{2} > 0, 
\end{equation}
with $z^i \in [-\varepsilon, \varepsilon]$.
If for example $\max(u(\cdot)-\bar{V}_\varepsilon (\cdot, 1))$ is reached in $x= -\varepsilon$ and $\max(u(\cdot)-\bar{V}_\varepsilon (\cdot, -1))$ is reached only in $\varepsilon$, then using \cref{uguaglianzavbarra} we get the contradiction
\begin{equation}
\begin{array}{ll}\label{staccodalbordo}
\displaystyle
u(-\varepsilon)- \bar{V}_{\varepsilon}(- \varepsilon, 1) = u(-\varepsilon)- \bar{V}_{\varepsilon}(- \varepsilon, -1)\\
\displaystyle
< u(\varepsilon)- \bar{V}_{\varepsilon}(\varepsilon, -1) = u(\varepsilon)- \bar{V}_{\varepsilon}(\varepsilon, 1).
\end{array}
\end{equation}
This implies that in at least one branch we can assume $z^i$ not equal to corresponding switching threshold. Then assume $z^{-1} \in [-\varepsilon, \varepsilon[$.

We are now comparing $u$ and $\overline V_\varepsilon$. By \cref{eq:vbarra}, we have for every $i=-1, 1$
\begin{equation*}
\lambda \bar{V}_{\varepsilon}(x, i)-f_{i}(x, \bar{a}_i)\bar{V}_{\varepsilon}^{'}(x, i)-\ell_i(x, \bar{a}_i)\geq - O(\varepsilon),
\end{equation*}
\noindent
in $x\in[-\varepsilon,\varepsilon[$ or in $]-\varepsilon,\varepsilon]$ respectively, where, here and in the sequel  $O(\varepsilon)$ is a suitable positive infinitesimal quantity as $\varepsilon \to 0$. Recalling that $\bar{V}_{\varepsilon}(\cdot, i)$ is derivable in $[-\varepsilon, \varepsilon]$ and recalling the sign of $f_i(0, \bar{a}_i)$ we then get for every $i=-1, 1$
\begin{equation}\label{Hamiltonianatrovata}
\lambda \bar{V}_{\varepsilon}(x, i)+ H_i(x, p)\geq - O(\varepsilon),
\end{equation}
for every $x \in [-\varepsilon, \varepsilon[$ and for every $p$ subgradient in $x$ with respect to $[-\varepsilon, \varepsilon]$ of $\bar{V}_{\varepsilon}(\cdot, -1)$ (respectively for any $ x \in ]-\varepsilon, \varepsilon]$ and $p$ subgradient of $\bar{V}_{\varepsilon}(\cdot, 1)$).

Let $\eta:[-\varepsilon, \varepsilon]\to \RR$ be continuous and $c>0$  such that (see \cite{So}, condition (A1) in the case of an interval)
\begin{equation} \label{funzioneeta}
 ]x+ \xi\eta(x)-\xi c, x+ \xi\eta(x)+\xi c [ \ \subseteq\ ]-\varepsilon, \varepsilon[ \ \forall x \in [-\varepsilon, \varepsilon], 0< \xi\leq c.
\end{equation}
For any $0<\xi\leq c$, we define the function in $[-\varepsilon, \varepsilon] \times [-\varepsilon, \varepsilon]$:
\begin{equation*}
\Phi_{\xi}(x, y) = u(x)- \bar{V}_{\varepsilon}(y, -1)- \left\vert\frac{x-y}{\xi}-\eta(z^{-1})\right\vert^2-\left\vert y-z^{-1}\right\vert^2. 
\end{equation*}
Let $(x_{\xi}^{-1}, y_{\xi}^{-1})$ be a point of maximum for  $\Phi_\xi \in [-\varepsilon, \varepsilon]\times [-\varepsilon, \varepsilon]$. Recalling $z^{-1} \in [-\varepsilon, \varepsilon[$ by standard estimates (see Soner \cite{So} or Bardi- Capuzzo Dolcetta \cite{BaCaDo} pp. 271) for small $\xi$ we get 
$x_\xi^{-1} \in ]-\varepsilon, \varepsilon[$, $y_\xi^{-1} \in [-\varepsilon, \varepsilon[$ and
\begin{equation}\label{stima2}
\frac{x_{\xi}^{-1} - y_{\xi}^{-1}}{\xi} \rightarrow \eta(z^{-1}) \quad \text{and}\quad x_{\xi}^{-1}, y_{\xi}^{-1}\rightarrow z^{-1} \quad\text{as}\quad \xi\rightarrow 0.
\end{equation}
\noindent We have the following possible cases, for a subsequence $\xi\to 0$:\\
(i) $(x_{\xi}^{-1}, y_{\xi}^{-1}) \in \ ]-\varepsilon,0[ \times [-\varepsilon, \varepsilon[$; (ii) $x_{\xi}^{-1}=0$ and $y_{\xi}^{-1}\in ]-\varepsilon, \varepsilon[$; (iii) $(x_{\xi}^{-1}, y_{\xi}^{-1}) \in \ ]0, \varepsilon[ \times ]-\varepsilon, \varepsilon[$.\\
Case (i). We get for any small $\xi$ 
\begin{equation}\label{sottosoluzine}
\lambda u(x_{\xi}^{-1})+ H_{-1} \bigg(x_{\xi}^{-1}, \frac{2}{\xi}\bigg(\frac{x_{\xi}^{-1} - y_{\xi}^{-1}}{\xi}-\eta(z^{-1})\bigg)\bigg)\leq 0,
\end{equation}
\begin{equation}\label{soprasoluzione}
\lambda \bar{V}_{\varepsilon}(y_{\xi}^{-1}, -1)+ H_{-1}\bigg(y_{\xi}^{-1}, \frac{2}{\xi}\bigg(\frac{x_{\xi}^{-1} - y_{\xi}^{-1}}{\xi}-\eta(z^{-1})\bigg)+2(z^{-1}-y_{\xi}^{-1})\bigg) \geq -O(\varepsilon).
\end{equation}
and we conclude in the standard way getting the contradiction to \cref{eq:ConditionMassimo} first sending $\xi \to 0$ and then $\varepsilon \to 0$.\\
Case (ii). By $x_{\xi}^{-1}=0$ we have that
\begin{equation}\label{eq:CondMinimo}
\begin{aligned}
\min\biggl\{\lambda u(0)+  H_1\bigg(0,\frac{2}{\xi}\bigg(\frac{- y_{\xi}^{-1}}{\xi}-& \eta(z^{-1})\bigg)\bigg),\\
&\lambda u(0)+ H_{-1}\bigg(0,\frac{2}{\xi}\bigg(\frac{- y_{\xi}^{-1}}{\xi}-\eta(z^{-1})\bigg)\bigg)\biggr\}\leq0.
\end{aligned}
\end{equation}
If $ \lambda u(0)+H_{-1}\bigg(0, \frac{2}{\xi}\bigg(\frac{- y_{\xi}^{-1}}{\xi}-\eta(z^{-1})\bigg)\bigg)\leq 0$ for a subsequence $\xi$ tends to 0 we conclude as in Case (i). Otherwise, we have
\begin{equation}\label{eq:CondSott}
\lambda u(0)+ H_1\bigg(0, \frac{2}{\xi}\bigg(\frac{- y_{\xi}^{-1}}{\xi}-\eta(z^{-1})\bigg)\bigg)\leq 0. 
\end{equation}
The inequality \cref{soprasoluzione,eq:CondSott} cannot be compared because they have different Hamiltonians. 
However, noting that $y_{\xi}^{-1} \in ]-\varepsilon, \varepsilon[$, we have 
\begin{equation*}
\left(\bar{V}_\varepsilon(y_{\xi}^{-1}, -1)\right)^{'}=\frac{2}{\xi}\bigg(\frac{x_{\xi}^{-1} - y_{\xi}^{-1}}{\xi}-\eta(z^{-1})\bigg)+2(z^{-1}-y_{\xi}^{-1}).
\end{equation*}

By \cref{DifferenzaNulla}, we have 
\begin{equation*}
\bar{V}_\varepsilon(y_{\xi}^{-1}, 1) = \bar{V}_\varepsilon(y_{\xi}^{-1}, -1)+ {O}(\varepsilon),\quad
 \left(\bar{V}_\varepsilon(y_{\xi}^{-1}, 1)\right)' = \left(\bar{V}_\varepsilon(y_{\xi}^{-1}, -1)\right)^{'}+ {O}(\varepsilon),
\end{equation*}
\noindent
and using \cref{Hamiltonianatrovata} in $y_{\xi}^{-1}$ for $i=1$,
we get
\begin{equation}\label{hamiltonianacontraria2}
\lambda \bar{V}_\varepsilon(y_{\xi}^{-1}, 1)+ H_1\bigg(y_{\xi}^{-1}, \frac{2}{\xi}\bigg(\frac{ - y_{\xi}^{-1}}{\xi}-\eta(z^{-1})\bigg)+2(z^{-1}-y_{\xi}^{-1}) \bigg)\geq - {O}(\varepsilon).
\end{equation}
By \cref{eq:CondSott,hamiltonianacontraria2} we obtain a contradiction as in the case (i).\\
 Case (iii). For $x_\xi^{-1} \in ]0, \varepsilon[$ we have
\begin{equation}\label{sottosoluzionemaggiorezero}
\lambda  u(x_{\xi}^{-1})+ H_{1} \bigg(x_{\xi}^{-1}, \frac{2}{\xi}\bigg(\frac{x_{\xi}^{-1}- y_{\xi}^{-1}}{\xi}-\eta(z^{-1})\bigg)\bigg)\leq 0
\end{equation}
that cannot be compared with \cref{soprasoluzione}. Being also $y_{\xi}^{-1} \in [-\varepsilon, \varepsilon[$, we conclude as before.
\end{proof}

\section{A threefold junctions problem}\label{sec:4}
\begin{figure}[htbp]
	\centering
	\includegraphics[scale=0.5]{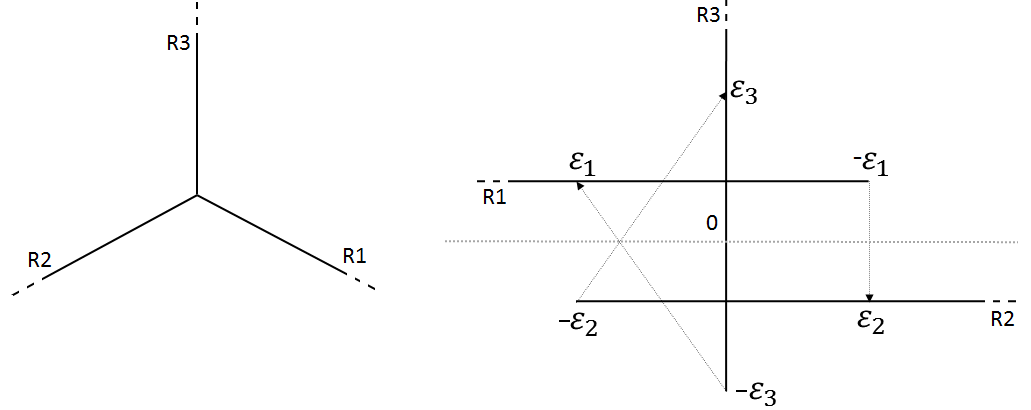}
\caption{\label{threejunction} The threefold junctions and its thermostatic-type approximation.}
\end{figure}

Here, we consider a junction given by three half-lines entering the same point (see \cref{threejunction}). In this case we have three labels $\left\{1, 2, 3\right\}$, one for every half-line $R_1, R_2, R_3$, that we identify with the labelled half-line  $R_i= [0, +\infty[\times \lbrace i \rbrace$. We also consider the controlled dynamics $f_i:R_i\times A \to \RR $ and the running costs $\ell_i: R_i\times A \to [0, +\infty[$. We approximate these triple discontinuity by a thermostatic-type combination in the following way. We  extend $f_i$ and $\ell_i$ to $[-\varepsilon_i, +\infty[\times \{i \} \times A$, where $\varepsilon_i>0$ are not necessarily the same for every $i$.
The thermostatic controlled dynamics is given by
\begin{equation}\label{s4}
\begin{cases}
x'(t) = f_{i(t)}(x(t), \alpha(t)),\\
i(t) = \tilde{h}[x](t), \\
i(0) = i_0 \in \left\{1, 2, 3\right\},\ 
x(0)=x_0 \in \ [-\varepsilon_{i_0}, +\infty[ ,
\end{cases}
\end{equation} 

\noindent where $\tilde{h}[x](t)$ is the delayed thermostatic rules as show \cref{threejunction}. In this thermostatic representation, denoting by $R_{\varepsilon_i}:=[-\varepsilon_i, +\infty[ \times \lbrace i \rbrace $ (and by $int(R_{\varepsilon_i})= ]-\varepsilon_i, +\infty[ \times \{i\}$), we can only switch from $R_{\varepsilon_1}$ to $R_{\varepsilon_2}$, from $R_{\varepsilon_2}$ to $R_{\varepsilon_3}$  and from $R_{\varepsilon_3}$ to $R_{\varepsilon_1}$. This is an arbitrary choice, because in the limit problem, at the junction-point,  a switch to any of the other branches is possible. However, we will recover this kind of behavior in the limit procedure because the transitions times become smaller and smaller and, indeed, the limit equation \cref{eq:HJBproblem3} is independent from that choice. Moreover, in the switching rule given by $\tilde{h}$, also the variable $x$ is subject to a discontinuity at the switching instant unlike the twofold case in previous section (see \cref{threejunction} and also note in the thermostat, the branch $R_{\varepsilon_1}$ is oriented in the opposite way with respect to the standard one). For every $i_0 \in \left\{1, 2, 3\right\}$ and $\forall x_0 \in \ [-\varepsilon_{i_0}, +\infty[$ we consider the value function 
\begin{equation}
V_{\varepsilon_1, \varepsilon_2, \varepsilon_3}(x_0, i_0)=\inf_{\alpha \in {\cal A}}\int_{0}^{\infty}e^{-\lambda t}\ell_{i(t)}(x(t), \alpha(t))dt,
\end{equation}
and we also have for every $i=1,2,3$ the Hamiltonians
\begin{equation}\label{s5}
H_i(x, p)= \sup_{a \in A}\left\{-f_i(x,a)\cdot p - \ell_i(x, a)\right\}.
\end{equation}
where we drop the index $i$ in the entries of $f_i$, $\ell_i$ and hence in $H_i$. We will sometimes use this simplification of the notation in the sequel too, without recalling it.\\ 
As in \cref{th:sistemapprox} we have the following proposition  
\begin{proposition}\label{Approx3}
For any choice of $\varepsilon_1, \varepsilon_2, \varepsilon_3>0$ the value function $V_{\varepsilon_1, \varepsilon_2, \varepsilon_3}$ of the switching three-thermostatic optimal control problem is the unique bounded and continuous function on $R_{\varepsilon_1} \cup R_{\varepsilon_2} \cup R_{\varepsilon_3}$ which satisfies, in the viscosity sense
\begin{equation}\label{RO}
\begin{cases}
\lambda V_{\varepsilon_1, \varepsilon_2, \varepsilon_3}(x, 1) + H_1 \bigg(x, V'_{\varepsilon_1, \varepsilon_2, \varepsilon_3}(x, 1)\bigg) = 0 & \text{in} \ \textit{int}(R_{\varepsilon_1}), \\
V_{\varepsilon_1, \varepsilon_2, \varepsilon_3}(-\varepsilon_1, 1) = V_{\varepsilon_1, \varepsilon_2, \varepsilon_3}(\varepsilon_2, 2); \\ 
\lambda V_{\varepsilon_1, \varepsilon_2, \varepsilon_3}(x, 2) + H_2 \bigg(x, V'_{\varepsilon_1, \varepsilon_2, \varepsilon_3}(x, 2)\bigg) = 0 & \text{in} \ \textit{int}(R_{\varepsilon_2}), \\
V_{\varepsilon_1, \varepsilon_2, \varepsilon_3}(-\varepsilon_2, 2) = V_{\varepsilon_1, \varepsilon_2, \varepsilon_3}(\varepsilon_3, 3); \\  
\lambda V_{\varepsilon_1, \varepsilon_2, \varepsilon_3}(x, 3) + H_3 \bigg(x, V'_{\varepsilon_1, \varepsilon_2, \varepsilon_3}(x, 3)\bigg) = 0 & \text{in} \ \textit{int}(R_{\varepsilon_3}),\\
V_{\varepsilon_1, \varepsilon_2, \varepsilon_3}(-\varepsilon_3, 3) = V_{\varepsilon_1, \varepsilon_2, \varepsilon_3}(\varepsilon_1, 1).
\end{cases}
\end{equation}
\end{proposition}
The proof is essentially an adaptation of the one of the \cref{th:sistemapprox} in \cite{Ba} to whom the reader is strongly referred. However, a very short sketch of the proof is given in the Appendix.

In the following subsections we are going to consider two different limit junction problems. In \cref{sectUni} the limit problem is given by the fact that on the junction point the admissible dynamics and costs are given by some suitably interpreted balance of the behaviours on all three branches. In \cref{subsec:nonuniform} instead, a balance among only two branches is also admitted. The main difference in the two limiting procedures is that in the first case the three thresholds go to zero with the same velocity, whereas in the second case different velocities are admitted. Such differences will lead to two HJB limit problems which differ by the definition of the admissible test functions (see the comments after \cref{def:visco}) which will lead to distinct ways for proving a comparison result.\\
	We finally note that a similar control problem with $n$ branches is studied in \cite{AcOuTc}. However, in that work no convexification or balance of the dynamics and costs are taken into account at the junction. The considered optimal control is then different from ours.	
\subsection{Uniform switching thresholds}\label{sectUni}
We assume $(\varepsilon_1, \varepsilon_2, \varepsilon_3)= (\varepsilon, \varepsilon, \varepsilon)$.
Looking to the twofold junction it is easy to see that the convexification parameters $\mu, 1-\mu$ are given by the ratio between the time spent using $f_i (0,a_i)$ to go from a threshold to the other one (namely $2\varepsilon/f_i (0,a_i) $) and the total time to perform a complete switching.
Coherently, when $f_1, f_2, f_3 < 0$ (dropping the entries in the dynamics), namely when we perform the whole cycle, the right convex parameters to be considered are 
\begin{equation}\label{musogliefisse}
\mu_1=\frac{f_{2}f_3}{f_{2}f_3 + f_{1}f_3 + f_{1}f_2}, \ \mu_2=\frac{f_1f_3}{f_{2}f_3 + f_{1}f_3 + f_{1}f_2}, \  \mu_3=\frac{f_1f_2}{f_{2}f_3 + f_{1}f_3 + f_{1}f_2}.
\end{equation}
Moreover $(\mu_1, \mu_2, \mu_3) \in [0, 1]^3$ and $\sum_{i=1}^{3} \mu_i = 1$.
Observe that now we have not anymore the interpretation as balance of forces, indeed in general $ \sum_{i=1}^{3} \mu_i f_i(0,a_i) \neq 0$, regardless to our choice of the signs of the branches $R_i$ and dynamics $f_i$. Also note that \cref{musogliefisse} is meaningful with the same interpretation when at most one $f_i$ is null, in which case we definitely remain in the corresponding branch.
To identify the right limit optimal control problem when $\varepsilon\to 0$ we define its controlled dynamics. In particular, calling $TR = R_1 \cup R_2 \cup R_3$, if $(x, i) \in TR$, with $x\neq 0$ then the dynamics is the usual $f_i(x, a_i)$ with $a_i\in A$. If instead $x=0$, being $(0, i)=(0,j)$ for $i, j \in \{1, 2, 3\}, i\neq j$,  we can either choose any dynamics makes us to stay inside a single branch $R_i$ or we may rest at zero ``formally" using any combination  $ \sum_{i=1}^{3} \mu_i f_i(0,a_i)$ with $f_i(0, a_i)$ and $\mu_i$ as before (where $\mu_i$ plays a role in the definition of the corresponding cost, see below). The set of controls in the junction point is then
\begin{equation*}
A(0)=A_{0}\cup \widetilde{A} 
\end{equation*}
with (note that in $\widetilde{A}$ the index $i$ is also at disposal)
\begin{align*}
A_{0} &= \lbrace (a_1, a_2, a_3) \in A^3 \vert \ f_i(0, a_i)\leq 0\ \text{with at most one equal to} \ 0\rbrace, \\
\widetilde{A} &= \left\{(a, i) \in A \times \{1, 2,3\} \vert\ f_i(0, a)\geq 0\ \right\}.
\end{align*}
Then, calling $\hat a$ the generic element of $A(0)$ we define 
\begin{equation*}
f_0(0, \hat a)=
\begin{cases}
f_i(0, a) & \text{if} \ \hat{a} \in  \widetilde{A}, \\
0 & \text{if}\ \hat{a} \in A_{0}.
\end{cases}
\end{equation*}
With the same arguments, if  $(x, i) \in TR$ and $x\neq 0$ then the running cost is  $\ell_i(x, a_i)$ with $a_i\in A$, otherwise we define
\begin{equation*}
\ell_0(0, \hat a)=
\begin{cases}
\ell_i(0, a) & \text{if} \ \hat{a} \in  \widetilde{A}, \\
\mu_1 \ell_1(0, a_1)+\mu_2\ell_2(0, a_2)+ \mu_3\ell_3(0, a_3) & \text{if}\ \hat{a} \in A_{0}.
\end{cases}
\end{equation*}
The quadruples $f = (f_1, f_2, f_3, f_0)$ and $\ell = (\ell_1, \ell_2, \ell_3, \ell_0)$ then define the threefold junction optimal control problem. In particular given an initial state $(x_0, i_0) \in TR$ and a measurable control $\alpha(t) \in A\cup A(0)$ we consider a possible admissible trajectory in $TR$ whose evolution, denoted by $(x(t), i(t))$, is such that $i(t)$ remains constant whenever $x(t)>0$ and $x(t)$ evolves with dynamics described above. Given an initial state, the set of measurable controls for which there exists a unique admissible trajectory is not empty and we denote it by ${\cal A}_{(x_0, i_0)}$.
We then consider an infinite horizon problem with a discount factor $\lambda >0$ given by
\begin{equation*}
J(x_0, i_0, \alpha) = \int_{0}^{+\infty} e^{-\lambda t}\ell(x(t), i(t), \alpha(t))  dt, 
\end{equation*}
where $\ell$ is the running cost described above and the corresponding
 value function is
\begin{equation}\label{funzionevaloresogliefisse}
V(x_0, i_0)= \inf_{\alpha \in {\cal A}_{(x_0, i_0)}}J(x_0, i_0, \alpha).
\end{equation}

In the sequel when $x=0$ we will drop the index $i$.
If we remain in $x=0$ for all the time using controls in $A_0$ then the best cost is given by
\begin{equation} \label{s15}
u_{1,2,3}(0) = \frac{1}{\lambda}\inf_{A_0}\left\{\mu_1\ell_1(0, a_1) + \mu_2\ell_2(0, a_2) + \mu_3\ell_3(0, a_3)\right\}.
\end{equation}
\begin{rem}\label{Remarkcontrolli}
In general $A_0$ is not compact. However, if $(a_1^{k}, a_2^{k}, a_3^{k}) \in A_0$ is a minimizing sequence for $u_{1,2,3}(0)$ converging to $(\bar{a}_1, \bar{a}_2, \bar{a}_3) \notin {A_0}$, the quantity inside the bracket in \cref{s15} loses meaning but we still have the inequality
\begin{equation*}
\lim_{k\rightarrow \infty}\left\{\mu_1^{k}\ell_1(0, a_1^{k}) + \mu_2^{k}\ell_2(0, a_2^k) + \mu_3^k \ell_3(0, a_3^k)\right\} \geq \min\{\ell_i(0, \bar{a}_i) \vert f_i(0,\bar{a}_i)=0\}. 
\end{equation*}
and we can still get an optimal behavior among the ones making us stay at $x=0$.
\end{rem}
\begin{theorem} \label{T1}
Assume \cref{eq:Lip,eq:Controllability,eq:LLip}. 
Then, $V$ is continuous on $TR$. Moreover when $x=0$,
\begin{equation} \label{s6}
V(0) = \min \left\{u_{1, 2, 3}(0), V_{sc(1)}(0), V_{sc(2)}(0), V_{sc(3)}(0) \right\},
\end{equation}
where $V_{sc(i)}(0)$ is the value function at $x=0$ of the state- constraint optimal control problem on $R_i$. Therefore\\
i) if $V(0) = u_{1, 2, 3}(0)$, then $V$ is the unique bounded and continuous solution of the three problems (one for every $i \in \left\{1, 2, 3\right\}$)
\begin{equation}
\begin{cases} \label{s7}
\lambda u + H_{i}(x, u^{\prime}) = 0 \ \ \text{in} \ int(R_i) \\
u(0) = u_{1, 2, 3}(0)
\end{cases}
\end{equation}
ii) if $V(0) = V_{sc(i)}(0)$, for some $i = 1, 2, 3$, then $V$ satisfies: $V = V_{sc(i)}$ in $R_i$, and uniquely solves (for every $j \in \left\{1, 2, 3\right\} \setminus \left\{i\right\}$)
\begin{equation}
\begin{cases} \label{s8}
\lambda u + H_{j}(x, u^{\prime}) = 0 \ \ \text{in} \  int(R_j) \\
u(0) = V_{sc(i)}(0).
\end{cases}
\end{equation}
\end{theorem}
\begin{proof}
The continuity of $V$ comes from controllability \cref{eq:Controllability} and regularity \cref{eq:Lip,eq:LLip} in a standard way.
Moreover, \cref{s6} comes from \cref{funzionevaloresogliefisse} because the four terms in the minimum are exactly the only allowed values (see also Remark \ref{Remarkcontrolli}). Finally \cref{s7,s8} follow from standard properties of Dirichlet problems in the viscosity sense. 
\end{proof}

\begin{theorem}\label{limitsogliefisse} Assume \cref{eq:Lip,eq:Controllability,eq:LLip}. The value function $V$ \cref{funzionevaloresogliefisse} (also characterized by \cref{T1}) satisfies 
\begin{equation} \label{s16}
V(x, i) = \lim_{\varepsilon \rightarrow 0}V_{\varepsilon, \varepsilon, \varepsilon}(x, i) \ \ \forall \ (x, i) \in  R_i, \ i= 1, 2, 3.
\end{equation}
where $V_{\varepsilon, \varepsilon, \varepsilon}$ is the value function of the approximating thermostatic problem \cref{RO} with uniform thresholds $(\varepsilon, \varepsilon, \varepsilon)$, and the convergence is uniform. Moreover, when $x=0$ the limit is independent from $i = 1, 2, 3$.
\end{theorem}
\begin{proof} 
We first prove that \cref{s16} holds for $x=0$ (the junction point).
The fact that the limit \cref{s16}, whenever it exists, is independent from $i$ when $x=0$ comes from the controllability hypothesis \cref{eq:Controllability} because $\vert V_{\varepsilon, \varepsilon, \varepsilon}(0, i) - V_{\varepsilon, \varepsilon, \varepsilon} (0, j) \vert$ is infinitesimal as $\varepsilon$. In the sequel, we drop the symbol $i$ in the expression $V_{\varepsilon, \varepsilon, \varepsilon}(0, i)$.\\
We prove \cref{s16} at $x=0$ for a convergent subsequence still denoted by $(\varepsilon, \varepsilon, \varepsilon)$ which exists because $ V_{\varepsilon, \varepsilon, \varepsilon}$ are equi-bounded. The uniqueness of the limit will give the whole \cref{s16}.
By contradiction, suppose that $V(0) < \lim V_{\varepsilon, \varepsilon, \varepsilon}(0)$. By \cref{eq:Controllability}, for every $\varepsilon> 0$, we have $V_{\varepsilon, \varepsilon, \varepsilon}(0) \leq V_{sc(i)}(0)$ for every $i=1, 2, 3$. Hence, the absurd hypothesis implies $V(0) = u_{1, 2, 3}(0)$ by \cref{s6}. 
Suppose that $(a_1,a_2, a_3) \in A_0$ realizes the minimum in the definition of $u_{1, 2, 3}(0)$.
We analyze some possible cases, the others are similar.

1) $f_1(0,a_1), f_2(0, a_2), f_3(0, a_3) < 0$. Hence, using a suitably switching control between those constant controls, we get  $V_{\varepsilon, \varepsilon, \varepsilon}(0) $ is not larger then $ u_{1, 2, 3}(0)$ plus an infinitesimal quantity as $\varepsilon\to 0$, which is a contradiction. 

2) $f_1(0,a_1) = 0, f_2(0, a_2), f_3(0, a_3) < 0$.
In this case we arrive at $R_1$ and we stop with $f_1(0,a_1)$ in $x=0$. Hence, $u_{1, 2, 3}(0)= \frac{1}{\lambda} \ell_1(0,a_1)$ cannot be lower than $V_{sc(1)}(0)$ which is a contradiction.

If there is no minimizing $(a_1,a_2,a_3)$ (see Remark \ref{Remarkcontrolli}) then the cost $u_{1, 2, 3}(0)$ cannot be better than the state-constraint value $V_{sc(i)}(0)$. Then as before, we have again a contradiction.

Now assume $\lim V_{\varepsilon, \varepsilon, \varepsilon}(0) < V(0)$. Let $\delta > 0$ be such that, for $\varepsilon$ small enough, it is $V_{\varepsilon, \varepsilon, \varepsilon}(0) + \delta < V(0)$. A measurable control $\alpha$ which almost realizes the optimum (less than $\beta>0$) for $V_{\varepsilon, \varepsilon, \varepsilon}(0)$ must be such that there are infinitely many switching between all branches $R_i^{\varepsilon}$ (i.e for every $i$, $f_i(x, \alpha_i)<0\ \forall\ x$ ).
Indeed, if it is not the case, then, for at least one branch $R_i^{\varepsilon}$, the trajectory definitely remains inside it. Hence, for small $\varepsilon$, $V_{\varepsilon, \varepsilon, \varepsilon}(0)$ is almost equal to $V_{sc(i)}(0)$, which is a contradiction.  We can restrict to consider a piecewise constant control that we call again $\alpha$ since $V_{\varepsilon, \varepsilon, \varepsilon}$ defined either by measurable controls or by piecewise constant controls, satisfies the same problem \cref{RO} which admits a unique solution (see e. g. \cite{BaCaDo} Remark 2.15 page 109). Then, to obtain the optimum, on each branch $R_i^{\varepsilon}$ let $x_1^i, \dots, x_{n^i}^i$ be the points corresponding to the discontinuity instants $t_1^i, \dots, t_{n^i}^i$ of the control $\alpha$ and let $a_j^i$ be the constant controls $\forall i=1, 2, 3$, $\forall \ j=1, \dots, n^i-1$. On the assumption that $f_i(0, a_j^i)<0 \ \forall\ i, j$ we consider the dynamics $f_i(0, a_j^i)$ and the running cost $\ell_i(0, a_j^i)$ on every spatial interval $[x_j^i, x_{j+1}^i]$. 
Now, for every $i$ we consider 
\begin{equation}\label{infdelrapporto}
\inf_{a \in A}\left\{ \frac{\ell_i(0, a)}{\vert f_i(0, a)\vert}\ \vert f_i(0, a)<0\right\}.
\end{equation}
If \cref{infdelrapporto} is a minimum for every $i$ obtained in $(\bar{a}_1, \bar{a}_2, \bar{a}_3)$ then in each $R_i^{\varepsilon}$ we use constant dynamics $f_i(0, \bar{a}_i)$ and constant running cost $\ell_i(0, \bar{a}_i)$. \\Therefore  $\left\vert J(\cdot, i, \alpha)- J(\cdot, i, \bar{a}_i)\right\vert \leq O(\varepsilon)$ and we get
\begin{equation}\label{primacontraddizione}
\begin{array}{ll}
\displaystyle
V_{\varepsilon, \varepsilon, \varepsilon}(0)\geq J(\cdot, i, \alpha)-\beta\geq J(\cdot, i, \bar{a}_i)-O(\varepsilon)-\beta\\
\displaystyle
\geq u_{1, 2, 3}(0)-O(\varepsilon)-\beta \geq V(0)-O(\varepsilon)-\beta,
\end{array}
\end{equation}
that is a contradiction. If, for some $i$, \cref{infdelrapporto} is not a minimum then we can consider the minimizing sequence $a_i^k$ that realizes the infimum less than $O(\frac{1}{k})$. In particular $a_i^k \to \tilde{a}_i \in A$ for $k\to +\infty$ and $f_i(0, a_i^k)\to f_i(0, \tilde{a}_i)=0$ being $f_i(0, a_i^k)<0$. However, since the optimal strategy is to switch among the branches, we cannot stop in the branch $R_i^{\varepsilon}$ with dynamics $f_i(0, \tilde{a}_i)$ paying the cost $\ell_i(0, \tilde{a}_i)$. Then, always taking into account that $f_i(0, a_j^i)<0$ we have
\begin{equation}\label{secondacontradd}
\begin{array}{ll}
\displaystyle
V_{\varepsilon, \varepsilon, \varepsilon}(0)\geq J(\cdot, i, \alpha)-\beta\geq J(\cdot, i,a_i^k)-O\left(\frac{1}{k}\right)-O(\varepsilon)-\beta\\
\displaystyle
\geq  u_{1, 2, 3}(0)-O\left(\frac{1}{k}\right)-O(\varepsilon)-\beta \geq V(0)-O\left(\frac{1}{k}\right)-O(\varepsilon)-\beta,
\end{array}
\end{equation}
which is again a contradiction.\\
Therefore at the end, $V_{\varepsilon, \varepsilon, \varepsilon}(0)$ cannot be less than $V(0) - \delta$ by the definition of $V(0)$. This is a contradiction. Hence we have $\lim V_{\varepsilon, \varepsilon, \varepsilon}(0) = V(0)$. The equations solved by $V_{\varepsilon,\varepsilon, \varepsilon}$ and by $V$ (\cref{RO,s7,s8} respectively) are the same for all $(x, i) \in int(R_i)$ and the boundary datum converges to $V(0)$. Hence, representing the solutions as the value functions of the corresponding optimal control problems, we get \cref{s16} and the uniform convergence.
\end{proof}
To show that $V$ \cref{funzionevaloresogliefisse} is a viscosity solution of the next problem \cref{eq:HJBproblem3}, we introduce the test functions for the differential equations on the branches and give the definition of viscosity subsolution and supersolution of \cref{eq:HJBproblem3}.

\begin{definition}\label{def_giromondo}
Let $\varphi: TR\to \RR$ be a function such that
\begin{equation}\label{test_function}
\begin{aligned}
&\varphi|_{R_i}:=\varphi_i : R_i \longrightarrow \RR \\
& (x, i) \longmapsto \varphi_i(x, i) \ \ \ \text{if}\ x \neq 0, \forall i \in \lbrace 1, 2, 3\rbrace \\
& (0, i) \longmapsto \varphi_i(0, i)=\varphi_j(0, j)\ \forall j \in \lbrace 1, 2, 3\rbrace \setminus \lbrace i \rbrace,
\end{aligned}
\end{equation}
with $\varphi \in C^0(TR)$ and $\varphi_i \in C^1(R_i)$.
\end{definition}
\begin{definition} \label{defsubsup}
A continuous function $u:TR\to \RR$ is a viscosity subsolution of \cref{eq:HJBproblem3} if for any $(x, i) \in TR$, any $\varphi$ as in \cref{test_function} such that $u-\varphi$ has a local maximum at $(x, i)$ with respect to $TR$, then 
\begin{equation}
\begin{aligned}
&\lambda u(x, i)+H_i(x, \varphi'_i(x, i))\leq 0 & & x \in int(R_i), \\
& \min\left\{\lambda u(0, i)+H_i(0, \varphi'_i(0, i)), \ i=1, 2, 3 \right\}\leq 0 & & x=0. 
\end{aligned}
\end{equation}
A continuous function $u:TR\to \RR$ is a viscosity supersolution of \cref{eq:HJBproblem3} if for any $(x, i) \in TR$, any $\varphi$ as in \cref{test_function} such that $u-\varphi$ has a local minimum at $(x, i)$ with respect to $TR$, then 
\begin{equation}
\begin{aligned}
&\lambda u(x, i)+H_i(x, \varphi'_i(x, i))\geq 0 & & x \in int(R_i), \\
& \max\left\{\lambda u(0, i)+H_i(0, \varphi'_i(0, i)),\ i=1, 2, 3 \right\}\geq 0 & & x=0.
\end{aligned}
\end{equation}
In particular note that if $x=0$ then the local maximum/minimum is with respect to all the three branches and $\varphi'_i(0, i)$ is the right derivative on the branch $i$, $(\varphi'_i)^+$. Since $(0, i) = (0, j)$ for $i,\ j \in \{1, 2, 3\} , \ i\neq j$, we drop the index $i$ in the pair $(0, i)$.
\end{definition}
We will prove the following theorem using the thermostatic approximation, namely considering the approximating value function $V_{\varepsilon, \varepsilon, \varepsilon}$. Differently from the twofold junction problem in which the index that identifies the branch is included in the sign of $x$ and the test function $\varphi \in C^1(\RR)$, here, we need to extend the test function $\varphi_i$ in \cref{test_function} from $R_i$ to $R_i^{\varepsilon}$. To do that we distinguish the case in which $V-\varphi$ has a local maximum point at $x=0$ from that where $x=0$ is a local minimum point, both respect to all three branches.\\
If $V-\varphi$ has a local maximum point at $x=0$ then we suppose that 
\begin{equation}\label{assumptionextendedtestfunction}
\varphi_1'(0)^+\leq\varphi_2'(0)^+\leq\varphi_3'(0)^+.
\end{equation}
Note that our switching sequence is $1\to 2\to 3\to 1$ which is coherent with such an order. If the order is different, then we consider a different switching sequence in the approximating thermostatic $\varepsilon$-problem, still coherent with the order. This is always possible because the limit function $V$ is independent from the switching order of the chosen approximating problem. Then we define
\begin{equation}\label{extendedtestfunction}
\tilde{\varphi}_i: [-\varepsilon, +\infty[\times \lbrace i\rbrace \longrightarrow \RR,\quad
\tilde{\varphi}_i =
\begin{cases} 
\varphi_i(x, i) &  x\geq 0 \\
\varphi_{i_s}(-x, i_s) &  x < 0
\end{cases}
\end{equation}
for $i=1, 2$ and with $i_s$ the next transition to $i$. If $i=3$ we construct $\tilde{\varphi}_3$ in two different way, correspondingly to the cases:
\begin{equation}\label{estensionevarphiramo3_1}
 \text{if} \quad \varphi_1'(0)^+=\varphi_3'(0)^+ \quad \text{then}
\quad \tilde{\varphi}_3=
\begin{cases}
\varphi_3(x, 3) &  x\geq 0, \\
\varphi_{1}(-x, 1) &  x < 0.
\end{cases}
\end{equation}
\begin{equation}\label{estensionevarphiramo3_2}
\text{if} \quad \varphi_1'(0)^+<\varphi_3'(0)^+ \quad \text{then}\quad
\tilde{\varphi}_3=
\begin{cases}
\varphi_3(x, 3) &  x\geq 0, \\
\varphi_{3}(-x, 3) &  x < 0.
\end{cases}
\end{equation}
By the assumption \cref{assumptionextendedtestfunction}, the first case gives, $\varphi_1'(0)^+=\varphi_2'(0)^+=\varphi_3'(0)^+$, and that the second case gives, at least for small $\varepsilon$, $\tilde\varphi_1(\varepsilon, 1)=\varphi_1(\varepsilon, 1)\le\varphi_3(\varepsilon, 3)=\tilde\varphi_3(-\varepsilon, 3)$. Finally in both cases we then have $\tilde\varphi_1(\varepsilon, 1)\le\tilde\varphi_3(-\varepsilon, 3)$.\\
If instead $V-\varphi$ has a local minimum point at $x=0$ then we suppose that 
\begin{equation}\label{assumptionextendedtestfunction2}
\varphi_1'(0)^+\geq\varphi_2'(0)^+\geq\varphi_3'(0)^+,
\end{equation}
and that the switching order is the coherent one, as above. In this case we construct $\tilde{\varphi}_3$ as in \cref{extendedtestfunction,estensionevarphiramo3_1,estensionevarphiramo3_2}(with the only difference of the case $\varphi_1'(0)^+<\varphi_3'(0)^+$ replaced by $\varphi_1'(0)^+>\varphi_3'(0)^+$). In this case, for at least small $\varepsilon$ it is $\tilde\varphi_1(\varepsilon, 1)\ge\tilde\varphi_3(-\varepsilon, 3)$.\\
The function $\tilde{\varphi}_i$ is not differentiable in $x=0$, hence we cannot write a unique HJB equation for the function $V_{\varepsilon, \varepsilon, \varepsilon}(\cdot, i)$ in branch $R_i^{\varepsilon}$. To overcome the problem of discontinuity of $\tilde{\varphi}_i'$ in $x=0$ we interpret the behaviour of the dynamic $f_i(x, a_i)<0$ for $x \in ]-\varepsilon, 0[$ as entering in the next branch of the switching rule. More precisely, considering for example the branches $R_1^{\varepsilon}$ and $R_2^{\varepsilon}$, we define the function ${V}_{\varepsilon, \varepsilon, \varepsilon}(x, 1)=: \widetilde{V}_{\varepsilon, \varepsilon, \varepsilon}(-x, 2)$, the dynamics $-{f}_1(x, a)=:\tilde{f}_2(-x, a)$ and the relative running costs ${\ell}_1(x, a)=:\tilde{\ell}_2(-x, a)$ for $x\in ]-\varepsilon, 0[$.
In this way, for any $x \in ]-\varepsilon, 0[$ a local maximum point of $V_{\varepsilon, \varepsilon, \varepsilon}(\cdot, 1)- \tilde{\varphi}_1(\cdot, 1)$, we get that $V_{\varepsilon, \varepsilon, \varepsilon}(\cdot, 1)$ satisfies
	\begin{equation*}
	\lambda \widetilde{V}_{\varepsilon, \varepsilon, \varepsilon}(-x, 2)+\sup_{a \in A}\left\{-\tilde{f}_2(-x, a)\varphi_2(-x, 2)'-\widetilde \ell_2(-x, a)\right\}\leq 0.
	\end{equation*}
	which is equivalent, for the considerations before, to
	\begin{equation*}
	\lambda V_{\varepsilon, \varepsilon, \varepsilon}(x, 1)+\sup_{a \in A}\left\{-f_1(x, a)\tilde{\varphi}_1(x, 1)'-\ell_1(x, a)\right\}\leq 0. 
	\end{equation*}
	The same ideas can be applied to the other pairs of branches $(R_2^{\varepsilon}, R_3^{\varepsilon})$ and $(R_3^{\varepsilon}, R_1^{\varepsilon})$.
\begin{theorem} \label{thm:ishiitre} Assume \cref{eq:Lip,eq:Controllability,eq:LLip}. The value function $V$ \cref{funzionevaloresogliefisse} is a viscosity solution of the Hamilton-Jacobi-Bellman problem
\begin{equation}\label{eq:HJBproblem3}
\begin{cases}
\lambda V + H_{1}(x, V') = 0 \ \ \text{in} \ int(R_1) \\
\lambda V + H_{2}(x, V') = 0 \ \ \text{in} \ int(R_2) \\
\lambda V + H_{3}(x, V') = 0 \ \ \text{in} \ int(R_3) \\
\min \left\{\lambda V + H_1, \lambda V + H_{2}, \lambda V + H_{3}\right\} \leq 0 \ \text{on} \ x=0 \\
\max \left\{\lambda V + H_1, \lambda V + H_{2}, \lambda V + H_{3}\right\} \geq 0 \ \text{on} \ x=0
\end{cases}
\end{equation}
\end{theorem}
\begin{proof}
From \cref{Approx3,limitsogliefisse} and by classical convergence result, we get the first three eqaution in \cref{eq:HJBproblem3}.\\
We now prove the fourth equation in \cref{eq:HJBproblem3}. Let $\varphi$ as given in \cref{test_function} such that $V-\varphi$ has a strict relative maximum at $x=0$ with respect all the three branches and consider the assumption \cref{assumptionextendedtestfunction}.
For every $i$ it is
\begin{equation}\label{condizionesottsoludynentranti}
\lambda V(0)+\sup_{a\in A, f_i(0,a)\ge0}\{-f_i(0,a)\varphi_i'^+(0)-\ell_i(0,a)\}\le 0.
\end{equation}
\noindent
Indeed, for every $\varepsilon>0$, and for every $t>0$, we have ($V_{\varepsilon,\varepsilon,\varepsilon}$ solves DPP; $\alpha$ p.c. stays for piece-wise constant control)
\begin{equation*}
V_{\varepsilon,\varepsilon,\varepsilon}(0,i)\le\inf_{\alpha\ p.c., f_i(0,\alpha)\ge0}\left(
\int_0^te^{-\lambda s}\ell_i(x(s),\alpha(s))ds+e^{-\lambda t}V_{\varepsilon,\varepsilon,\varepsilon}(x(t),i)
\right)
\end{equation*}
\noindent
and hence, passing to the limit $\varepsilon\to0^+$,
\begin{equation*}
V(0)\le\inf_{\alpha\ p.c., f_i(0,\alpha)\ge0}\left(
\int_0^te^{-\lambda s}\ell_i(x(s),\alpha(s))ds+e^{-\lambda t}V(x(t),i)
\right)
\end{equation*}
\noindent
and finally we get the desired inequality \cref{condizionesottsoludynentranti}, being $x=0$ a local maximum for $V-\varphi_i$ with respect to $R_i$.
 Hence, we only need to prove that, with our hypotheses, for at least one $i$, we get
\begin{equation}\label{disequazionedaprovare}
\lambda V(0)+\sup_{a\in A, f_i(0,a)\le0}\{-f_i(0,a)\varphi_i'(0)^+-\ell_i(0,a)\}\leq 0.
\end{equation}
For each $i$ let  $(x_\varepsilon^i, i)$ be a sequence of local maximum points for $V_{\varepsilon,\varepsilon,\varepsilon}-\tilde\varphi_i$ with respect to $R_i^\varepsilon$ convergent to $(0, i)$, with $\tilde\varphi_i$ as in \cref{extendedtestfunction}. For each $\varepsilon$, for at least one branch $i$ we may assume $x_\varepsilon^i\neq-\varepsilon$. Indeed, if it is not the case, recalling that, by controllability implies $V_{\varepsilon,\varepsilon,\varepsilon}(-\varepsilon,i)\le V_{\varepsilon,\varepsilon,\varepsilon}(\varepsilon,i_s)$, we get the contradiction
\begin{equation*}
\begin{array}{ll}
V_{\varepsilon,\varepsilon,\varepsilon}(\varepsilon,1)-\tilde\varphi_1(\varepsilon, 1)<V_{\varepsilon,\varepsilon,\varepsilon}(-\varepsilon,1)-\tilde\varphi_1(-\varepsilon, 1)\le
V_{\varepsilon,\varepsilon,\varepsilon}(\varepsilon,2)-\tilde\varphi_2(\varepsilon, 2)< \\ V_{\varepsilon,\varepsilon,\varepsilon}(-\varepsilon,2)-\tilde\varphi_2(-\varepsilon, 2)\le
V_{\varepsilon,\varepsilon,\varepsilon}(\varepsilon,3)-\tilde\varphi_3(\varepsilon, 3)<V_{\varepsilon,\varepsilon,\varepsilon}(-\varepsilon,3)-\tilde\varphi_3(-\varepsilon, 3)\le\\
V_{\varepsilon,\varepsilon,\varepsilon}(\varepsilon,1)-\tilde\varphi_1(\varepsilon, 1).
\end{array}
\end{equation*}

Now, let $i$ be such that $x_\varepsilon^i\neq-\varepsilon$ for every $\varepsilon$ (or at least for a subsequence). If $x_\varepsilon^i>0$ for all $\varepsilon$, in the limit we get
\begin{equation*}
\lambda V(0)+\sup_{a\in A}\left\{-f_i(0,a)\varphi_i'(0)^+-\ell_i(0,a)\right\}\leq 0,
\end{equation*}
\noindent
and we get the conclusion.
\noindent
If $x_\varepsilon^i\in]-\varepsilon,0[$, in the limit we get

\begin{equation*}
\lambda V(0)+\sup_{a\in A}\left\{-f_i(0,a)\tilde\varphi_i'(0)^--\ell_i(0,a)\right\}\leq 0,
\end{equation*}

\noindent
and in particular
\begin{equation*}
\lambda V(0)+\sup_{a\in A,f_i(0,a)\le0}\left\{-f_i(0,a)\tilde\varphi_i'(0)^--\ell_i(0,a)\right\}\leq 0.
\end{equation*}

\noindent
where $\tilde\varphi_i'(0)^-$ is the left derivative of $\tilde\varphi_i$ at $x=0$. Now, if we are in the first case (all the right derivatives coincide) then we have $\tilde\varphi_i'(0)^-=\varphi_{i_s}'(0)^+=\varphi_i'(0)^+$, and hence we get \cref{disequazionedaprovare}. If instead $\varphi_1'(0)^+<\varphi_3'(0)^+$ then if $i=1$ then $i_s=2$, hence, by our hypotheses, in the inequality above it is 
\begin{equation*}
-f_i(0,a)\tilde\varphi_i'(0)^-=-f_i(0,a)\varphi_{i_s}'(0)^+\ge-f_i(0,a)\varphi_i'(0)^+,
\end{equation*}
\noindent
and we conclude. Same arguments if $i=2$ and $i_s=3$. If instead $i=3$, then $\tilde\varphi_3(0)^-=\varphi_3'(0)^+$ and we conclude.

Finally, if $x_\varepsilon^i=0$, then we still get 
\begin{equation*}
\lambda V(0)+\sup_{a\in A,f_i(0,a)\le0}\left\{-f_i(0,a)\tilde\varphi_i'(0)^--\ell_i(0,a)\right\}\le0,
\end{equation*}
\noindent
and we conclude as before, i.e. studying the two cases as above.\\
Now we suppose $V-\varphi$ have a local minimum with respect to $TR$ at $(0, i)$ and consider \cref{assumptionextendedtestfunction2}. We have to prove that, for at least one $i$, we have
\begin{equation}\label{disequazionedaprovare2}
\lambda V(0)+\sup_{a\in A}\{-f_i(0,a)\varphi_i'(0)^+-\ell_i(0,a)\}\ge0.
\end{equation}

If for some $i$ and for $\varepsilon\to0^+$, $V_{\varepsilon, \varepsilon, \varepsilon}(\cdot,i)$ coincides with the state-constraint value function on $R_\varepsilon^i$, then $V_{\varepsilon, \varepsilon, \varepsilon}(\cdot,i)$ and $V(\cdot, i)$ coincides on $R_i$ and hence $V$ satisfies the same HJB equation as $V_{\varepsilon, \varepsilon, \varepsilon}$, which is \cref{disequazionedaprovare2}.
\noindent
Hence we suppose that no $V_{\varepsilon, \varepsilon, \varepsilon}(\cdot,i)$ coincide with the corresponding state-constraint value function. 

For each $i$ let $(x_\varepsilon, i)$ be a sequence of local minimum points for $V_{\varepsilon, \varepsilon, \varepsilon}-\tilde\varphi_i$ with respect to $R_\varepsilon^i$, which converges to $(0, i)$. In this case we may assume that, for a fixed $i$, the sequence is such that either $x_\varepsilon^i\neq-\varepsilon$ or $x_\varepsilon^i=-\varepsilon$ but the HJB equation satisfied by $V_{\varepsilon, \varepsilon, \varepsilon}$ has the right sign ($\geq 0$). Indeed, if it is not the case (i.e. $x_\varepsilon^i=-\varepsilon$ and HJB has the wrong sign), we must have $V_{\varepsilon, \varepsilon, \varepsilon}(-\varepsilon,i)=V_{\varepsilon, \varepsilon, \varepsilon}(\varepsilon,i_s)$, and hence we get the following contradiction
\[
\begin{array}{ll}
V_{\varepsilon, \varepsilon, \varepsilon}(\varepsilon, 1)-\tilde\varphi_1(\varepsilon, 1)>V_{\varepsilon, \varepsilon, \varepsilon}(-\varepsilon,1)-\tilde\varphi_1(-\varepsilon, 1)=
V_{\varepsilon, \varepsilon, \varepsilon}(\varepsilon,2)-\tilde\varphi_2(\varepsilon, 2)>\\ V_{\varepsilon, \varepsilon, \varepsilon}(-\varepsilon, 2)-\tilde\varphi_2(-\varepsilon, 2)=
V_{\varepsilon, \varepsilon, \varepsilon}(\varepsilon, 3)-\tilde\varphi_3(\varepsilon, 3)>V_{\varepsilon, \varepsilon, \varepsilon}(-\varepsilon, 3)-\tilde\varphi_3(-\varepsilon, 3)\ge\\
V_{\varepsilon, \varepsilon, \varepsilon}(\varepsilon, 1)-\tilde\varphi_1(\varepsilon, 1).
\end{array}
\]
If $x_\varepsilon^i>0$, in the limit we exactly get
\[
\lambda V(0)+\sup_{a\in A}\{-f_i(0,a)\varphi_i'(0)^+-\ell_i(0,a)\}\ge0.
\]
If $x_\varepsilon^i\in[-\varepsilon,0[$ in the limit we get
\[
\lambda V(0)+\sup_{a\in A}\{-f_i(0,a)\tilde\varphi_i'(0)^--\ell_i(0,a)\}\ge0.
\]
If all the right derivatives at $x=0$ of $\varphi_i(0)$ coincide, then we conclude because $\tilde\varphi_i'(0)^-=\varphi_{i_s}'(0)^+=\varphi_i'(0)^+$ . Otherwise, if $i=3$ then we have $\tilde\varphi_i'(0)^-=\varphi_3'(0)^+$ and we conclude; if $i=1$ or $i=2$, by the hypotheses on $V_{\varepsilon, \varepsilon, \varepsilon}(\cdot,i)$ not coincident with the state-constraint value function, we get that the supremum above is approximated by controls such that $f_i(0,a)\le0$, which means 
\begin{equation*}
-f_i(0,a)\tilde\varphi_i'(0)^-=-f_i(0,a)\tilde\varphi_{i_s}(0)^+\le-f_i(0,a)\varphi_i(0)^+
\end{equation*}
\noindent
and we conclude.

If $x_\varepsilon^i=0$, then in the limit we get (still recalling that $V_{\varepsilon, \varepsilon, \varepsilon}(\cdot,i)$ is not the state-constraint value function)
\begin{equation*}
\lambda V(0)+\sup_{a\in A,f_i(0,a)\le0}\{-f_i(0,a)\tilde\varphi_i'(0)^--\ell_i(0,a)\}\ge0
\end{equation*}
\noindent
and we conclude as before.
\end{proof}
Now we want to prove that $V$ \cref{funzionevaloresogliefisse} is the maximal subsolution of \cref{eq:HJBproblem3}.
Assume that $\forall \ \varepsilon >0$ small enough, the optimal strategy for the approximating problem $\varepsilon$, starting by any $(x, i)$ with $x \in [-\varepsilon, \varepsilon]$, is  to run through infinitely many switches between the three branches (i.e. no state-constraint behaviour is optimal). Let $\mu_1, \mu_2, \mu_3$ be as in \cref{musogliefisse} and $(a_1, a_2, a_3) \in A_0$ realize the minimum in \cref{s15} such that
\begin{equation}\label{eq:Optimalitycaso3}
V(0)=u_{1,2,3}(0)= \frac{1}{\lambda}\lbrace \mu_1 \ell_1(0, a_1) + \mu_2 \ell_2(0, a_2)+ \mu_3 \ell_3(0, a_3) \rbrace.
\end{equation}
For every $x \in [0, \varepsilon]$, we define the following functions 
\begin{equation}\label{vbarracasotre}
\bar{V}^{\varepsilon}(x, i) = \int_{0}^{\frac{x}{\vert f_i(0, a_i)\vert}} e^{-\lambda t}\ell_i(0, a_i)dt + e^{\frac{-\lambda x}{\vert f_i(0, a_i)\vert}} u_{1, 2, 3}(0),
\end{equation}
where the upper extremal of the integration is the reaching time of the point $0$ in the corresponding branch starting from $x \in [0, \varepsilon]$. For $x \in [0, \varepsilon]$,  $V_{\varepsilon, \varepsilon, \varepsilon,}(x, i)$ is not larger then $\bar{V}^{\varepsilon}(x, i)$ plus an infinitesimal quantity as $\varepsilon\to 0$. The functions in \cref{vbarracasotre} are differentiable in $[0, \varepsilon]$. A direct computation gives
\begin{equation*}
\bar{V}^{\varepsilon}(x, i)' = \frac{\ell_i(0, a_i)}{\vert f_i(0, a_i)\vert}e^{\frac{-\lambda x}{\vert f_i(0, a_i)\vert}}-\frac{\lambda e^{\frac{-\lambda x}{\vert f_i(0, a_i)\vert}}}{\vert f_i(0, a_i)\vert}{u}_{1, 2, 3}(0),
\end{equation*}
and then for $\varepsilon\to 0$
\begin{align*}
\bar{V}^{\varepsilon}(x, 1)' & \longrightarrow \frac{(1-\mu_1)\ell_1(0, a_1)-\mu_2\ell_2(0, a_2)-\mu_3\ell_3(0, a_3)}{\vert f_1(0, a_1)\vert}, \\
\bar{V}^{\varepsilon}(x, 2)' & \longrightarrow \frac{-\mu_1\ell_1(0, a_1)+(1-\mu_2)\ell_2(0, a_2)-\mu_3\ell_3(0, a_3)}{\vert f_2(0, a_2)\vert}, \\
\bar{V}^{\varepsilon}(x, 3)' & \longrightarrow \frac{-\mu_1\ell_1(0, a_1)-\mu_2\ell_2(0, a_2)+(1-\mu_3)\ell_3(0, a_3)}{\vert f_3(0, a_3)\vert}.
\end{align*}
Moreover by \cref{vbarracasotre} we have for every $i=1, 2, 3$
\begin{equation}\label{condsoprasoluzionevbarra}
\lambda \bar{V}^{\varepsilon}(x, i)-f_i(x, a_i)\bar{V}^{\varepsilon}(x, i)'-\ell_i(x, a_i)\geq -O(\varepsilon),
\end{equation}
in $x \in [0, \varepsilon]$. In \cref{condsoprasoluzionevbarra} when $x=0$ we use the right derivative of $\bar{V}^{\varepsilon}(x, i)$  and $\bar{V}^{\varepsilon}(0, i)= {u}_{1, 2, 3}(0)$ for every $i$. Furthermore, by differentiability of $\bar{V}^{\varepsilon}(x, i)$ and recalling the sign of $f_i(0, a_i)$ we then get for every $i$
\begin{equation*}
\lambda \bar{V}^{\varepsilon}(x, i)+H_i(x, q)\geq -O(\varepsilon),
\end{equation*}
for every $x \in [0, \varepsilon]$ and for every $q$ subgradient in $x$ of $\bar{V}^{\varepsilon}(x, i)$.\\
We now define on $TR\cap \left(\cup_{i=1}^3 [0, \varepsilon] \times \{i\} \right)$ the function
\begin{equation}\label{vbarracasotrederivabile}
\bar{V}(x)=
\begin{cases}
\bar{V}^{\varepsilon}(x, i) &  \text{if} \ x \in int(R_i),\\
u_{1, 2, 3}(x) &  \text{if} \ x =0.
\end{cases}
\end{equation}
which is in $C^1([0, \varepsilon])$ and that we extend to whole $TR$ maintaining its differentiability. 
\begin{theorem} For each $u$ bounded, continuous subsolution of \cref{eq:HJBproblem3}, it is $u\leq V$ in $TR$.
\end{theorem}
\begin{proof} 
We can assume to be in the settings above for which \cref{eq:Optimalitycaso3} holds.  Indeed, otherwise in at least one branch $V$ coincides with the corresponding state-constraint value function, greater than any subsolution (see Soner \cite{So}). We then also have $u\leq V$ on the other branches. By contradiction, suppose $\sup_{(x, i) \in TR} (u-V)(x, i) > \delta >0$. If 
\begin{equation*} 
\exists r>0 \vert \forall \delta' >0\ \exists\ (\overline{x}, i) \in ]r, +\infty[\times \{i\}: \sup_{(x, i)\in TR} \big((u-{V})(x, i)-(u-V)(\overline{x}, i)\big)\leq \delta',
\end{equation*}
then, by
\cref{thm:ishiitre} and known comparison techniques we get a contradiction because, in $ ]r, +\infty[\times\{i \}$, $V$ is a supersolution and $u$ is a subsolution of the same HJB. Hence we may restrict to the case where $u -V$ has the maximum with respect to $r$ in $x=0$.
Since $\bar{V}^{\varepsilon}(x, i)$ converges to $V(0)$, with $\bar{V}^{\varepsilon}$ defined in \cref{vbarracasotre}, then for small $\varepsilon$,
\begin{equation}\label{condmassimo3}
u(z^i, i)-\bar{V}^{\varepsilon}(z^i, i)=\max_{[0, \varepsilon]\times\{i \}}(u(\cdot, i)-\bar{V}^{\varepsilon}(\cdot, i))>\frac{\delta}{2}>0.
\end{equation}
Since $u(x, i)$ is a continuous subsolution of \cref{eq:HJBproblem3} then satisfies
\begin{equation}
\lambda u(x, i)-f_i(x, a_i)\cdot p-\ell_i(x, a_i)\leq 0 \quad \forall p \in D^{+}u(x, i)\neq \emptyset,
\end{equation}
where $D^{+}u(x, i)$ is the set of super-differentials of $u$ at a point $(x, i)$.
Now, taking into account \cref{condsoprasoluzionevbarra,condmassimo3} we have that
\begin{equation}
p-\bar{V}^{\varepsilon}(x, i)' \leq \frac{-\lambda \delta}{2\vert f_i(x, a_i)\vert} +O(\varepsilon),
\end{equation}
whence, for $\varepsilon < \frac{1}{2}\left\vert \frac{\lambda \delta}{2\vert f_i(x, a_i)\vert}\right\vert$, we get that $p-\bar{V}^{\varepsilon}(x, i)'\leq -\bar{\delta}$, for a suitable $\bar{\delta}>0$ regardless to $x$. Hence $u(x, i)-\bar{V}^{\varepsilon}(x, i)$ is decreasing and, taking $\varepsilon$ as above, has maximum point in $x=0$.
By the previous consideration we get that $\bar{V}(x)$ \cref{vbarracasotrederivabile} is an admissible test function and that $u-\bar{V}$ has a local maximum point in $x=0$ for
suitable small $\varepsilon>0$. 
Hence, being $u$ a subsolution, exists $\bar{i} \in \{1, 2, 3\}$ such that
\begin{equation}\label{sottosoluzioneibarra}
\lambda u(0)+H_{\bar{i}}\left(0, \left(\bar{V}^{\varepsilon}(0, \bar{i})'\right)\right)\leq 0.
\end{equation}
Moreover, by\cref{vbarracasotre}, we have
\begin{equation}\label{soprasoluzioneibarra}
\lambda \bar{V}^{\varepsilon}(0, \bar{i})+H_{\bar{i}}\left(0, \left(\bar{V}^{\varepsilon}(0, \bar{i})'\right)\right)\geq -O(\varepsilon).
\end{equation}
Subtracting \cref{soprasoluzioneibarra} to \cref{sottosoluzioneibarra} we contradict \cref{condmassimo3} and then, for $\varepsilon\to 0$, $u\leq V$ in $TR$.
\end{proof}
\subsection{Non-uniform switching thresholds}
\label{subsec:nonuniform}
In this section we suppose that the three thresholds of the three-thermostatic optimal control problem are not the same for all $R_{\varepsilon_{i}}$. This imply that the time spent in a single branch $R_{\varepsilon_{i}}$ to reach the relative threshold depends on the value of $\varepsilon_i$. Accordingly to this, the convexification parameters $\bar{\mu}_1, \bar{\mu}_2, \bar{\mu}_3$ are such that if at limit for $(\varepsilon_1, \varepsilon_2, \varepsilon_3)\to (0^+, 0^+, 0^+)$ the optimal behavior is to switch only between two branch, $R_i$ and $R_j$ for $i, j\in \{1, 2, 3\}, i\neq j$, then $\bar{\mu}_i+\bar{\mu}_j=1$. If instead the optimal  behavior is to switch among all three branches $R_i$ then $\bar{\mu}_i=\mu_i$ as in \cref{musogliefisse}. To identify the limit optimal control problem when $(\varepsilon_1, \varepsilon_2, \varepsilon_3) \to (0^+, 0^+, 0^+)$ we define the controlled dynamics. Using the same notation of the last section, if $(x, i) \in TR$ with $x\neq 0$ then the dynamics is the usual $f_i(x, a_i)$ with $a_i \in A$. If instead $x=0$, being $(0, i)=(0,j)$ for $i, j \in \{1, 2, 3\}, i\neq j$,  we can either choose any dynamics makes us to stay inside a single branch $R_i$ or we may rest at zero using any combination  $ \sum_{i=1}^{3} \bar{\mu}_i f_i(0,a_i)$ with $f_i(0, a_i)$ and $\bar{\mu}_i$ as before. In detail,  the set of controls in the junction is $A(0)=\overline{A}\cup \widetilde{A}$ 
with
\begin{align*}
\overline{A} & = \{(a_1, a_2, a_3, \sigma, \bar{\mu}_1, \bar{\mu}_2, \bar{\mu}_3 ) \in A^3 \times \{12, 13, 23, 123 \} \times [0, 1]^3 \vert \\
& \sigma=ij \Rightarrow \bar{\mu}_i+\bar{\mu}_j=1, f_i(0, a_i)\leq 0;\\
& \sigma=123 \Rightarrow \bar{\mu}_i=\mu_i, f_i(0, a_i)\leq 0 \ \text{with at most one equal to} \ 0 \}, \\
\widetilde{A} &= \left\{(a, i) \in A \times \{1, 2,3\} \vert\ f_i(0, a)\geq 0\right\}.
\end{align*}
In $\widetilde{A}$ the index $i$ is at disposal, while in $\overline{A}$, the notation $ij$ means that the switching is only between $R_i$ and $R_j$ (as well as $123$ means that the switching performs among all the three branches).\\
Then, as in the last section, calling $\hat{a}$ the generic element of $A(0)$ we define
\begin{equation*}
f_0(0, \hat a)=
\begin{cases}
f_i(0, a) & \text{if} \ \hat{a} \in  \widetilde{A}, \\
0 & \text{if}\ \hat{a} \in \overline{A}.
\end{cases}
\end{equation*}
With the same arguments, if  $(x, i) \in TR$ and $x\neq 0$ then the running cost is  $\ell_i(x, a_i)$ with $a_i\in A$, otherwise we define
\begin{equation*}
\ell_0(0, \hat a)=
\begin{cases}
\ell_i(0, a) & \text{if} \ \hat{a} \in  \widetilde{A}, \\
\bar{\mu}_1 \ell_1(0, a_1) + \bar{\mu}_2 \ell_2(0, a_2) & \text{if} \ \sigma=12 \quad \text{and}\ \hat{a} \in \overline{A},\\
\bar{\mu}_1 \ell_1(0, a_1) + \bar{\mu}_3 \ell_3(0, a_3) & \text{if} \ \sigma=13 \quad \text{and}\ \hat{a} \in \overline{A}, \\
\bar{\mu}_2 \ell_2(0, a_2) + \bar{\mu}_3 \ell_3(0, a_3) & \text{if} \ \sigma=23 \quad \text{and}\ \hat{a} \in \overline{A} \\
\mu_1 \ell_1(0, a_1)+\mu_2\ell_2(0, a_2)+ \mu_3\ell_3(0, a_3) & \text{if} \ \sigma=123 \quad \text{and}\ \hat{a} \in \overline{A}.
\end{cases}
\end{equation*}


The quadruples $f = (f_1, f_2, f_3, f_0)$ and $\ell = (\ell_1, \ell_2, \ell_3, \ell_0)$ then define a threefold junction optimal control problem, different from the one in \cref{subsec:nonuniform}. We still denote by ${\cal A}_{(x_0, i_0)}$ the nonempty set of measurable controls for which there exists a unique admissible trajectory and consider the cost functional ($\lambda >0$)
\begin{equation*}
J(x_0, i_0, \alpha) = \int_{0}^{+\infty} e^{-\lambda t} \ell(x(t), i(t), \alpha(t))  dt
\end{equation*}
where $\ell$ is the running cost described above. The corresponding value function is
\begin{equation}\label{funzionevaloresogliemobili}
V^*(x_0, i_0)= \inf_{\alpha \in {\cal A}_{(x_0, i_0)}}J(x_0, i_0, \alpha).
\end{equation}
Observe that if we stay in $x=0$ for all time using controls in $\overline{A}$ the cost is 
\begin{equation*}
u_{0}(0) = \frac{1}{\lambda}\ \min_{\overline A}\ \sum_{i=1}^{3} \bar{\mu}_i \ell_{i}(0, a_i) = \frac{1}{\lambda}\ \min\left\{u_{1, 2}(0), u_{1, 3}(0), u_{2, 3}(0), u_{1, 2, 3}(0)\right\}
\end{equation*}
where $u_{1, 2}(0)$ is the minimum over ${\overline A}$  of the cost $\ell_0$ when $\sigma=12$, $u_{1, 3}(0)$  is the minimum over ${\overline A}$ of the cost $\ell_0$ when $\sigma=13$ and similarly the others.
\begin{theorem}\label{limitesogliemobili} Assume \cref{eq:Lip,eq:Controllability,eq:LLip}. The value function $V^*$ \cref{funzionevaloresogliemobili} characterized by \cref{T1}, but with $u_0(0)$ in place of $u_{1, 2, 3}(0)$, namely
\begin{equation} \label{defVstarinzero}
V^*(0) = \min \left\{u_{0}(0), V_{sc(1)}(0), V_{sc(2)}(0), V_{sc(3)}(0) \right\},
\end{equation}
satisfies
\begin{equation} \label{s11}
V^{*}(x, i) = \liminf_{(\varepsilon_1, \varepsilon_2, \varepsilon_3) \rightarrow (0^{+}, 0^{+}, 0^{+})} V_{\varepsilon_1, \varepsilon_2, \varepsilon_3}(x, i) \ \forall \ (x, i) \in  R_i, \ i= 1, 2, 3.
\end{equation}
where $V_{\varepsilon_1, \varepsilon_2, \varepsilon_3}$ is the value function of the approximating thermostatic problem \cref{RO}, with non uniform thresholds $(\varepsilon_1, \varepsilon_2, \varepsilon_3)$, and the convergence is uniform. Moreover, when $x=0$, the limit is independent from $i=1, 2, 3$.
\end{theorem}
\begin{proof}
We prove \cref{s11} at $x=0$. The independence from $i$ of  \cref{s11} comes from the controllability \cref{eq:Controllability}: $|V_{\varepsilon_1, \varepsilon_2, \varepsilon_3}(0, i) - V_{\varepsilon_1, \varepsilon_2, \varepsilon_3}(0, j)|$\ is infinitesimal as $\max \left\{\varepsilon_1, \varepsilon_2, \varepsilon_3\right\}$. In the sequel, we omit the symbol $i$ in the expression $V_{\varepsilon_1, \varepsilon_2, \varepsilon_3}(0, i)$.\\
By contradiction, suppose $V^{*}(0) < \liminf V_{\varepsilon_1, \varepsilon_2, \varepsilon_3}(0)$. By \cref{eq:Controllability}, for every $\varepsilon_1, \varepsilon_2, \varepsilon_3 > 0$, we have $V_{\varepsilon_1, \varepsilon_2, \varepsilon_3}(0) \leq V_{sc(i)}(0)$ for every $i=1, 2, 3.$ Hence, it implies $V^{*}(0) = u_{0}(0)$. Let $(a_1, a_2, a_3, \sigma, \bar{\mu}_1, \bar{\mu}_2, \bar{\mu}_3) \in \overline{A}$ realize the minimum in the definition of $u_{0}(0)$. We analyze some possible cases, the other ones being similar.

1) $f_1(0,a_1), f_2(0, a_2), f_3(0, a_3) < 0$ and $\sigma=123$. Taking $(\varepsilon_1, \varepsilon_2, \varepsilon_3) = (\varepsilon, \varepsilon, \varepsilon)$ and using a suitably switching control between those constants controls, $ V_{\varepsilon, \varepsilon, \varepsilon}(0)$ is not larger than $u_{1, 2, 3}(0)$ plus an infinitesimal quantity as $\varepsilon\to 0$, which is a contradiction.

2) $f_1(0,a_1), f_2(0, a_2), f_3(0, a_3) < 0$ and $\sigma=23$. Here, taking the triple \\ $(\varepsilon_1, \varepsilon_2, \varepsilon_3) = (\varepsilon^2, \varepsilon, \varepsilon)$, $V_{\varepsilon, \varepsilon, \varepsilon}(0)$ is not larger than $u_{2, 3}(0)$ plus an infinitesimal quantity as $\varepsilon\to 0$, which is a contradiction. 

3) $f_1(0,a_1) = 0, f_2(0, a_2), f_3(0, a_3) < 0$. In this setting we can study two sub-cases according to the value of $\sigma$.

3.1) If $\sigma=123$, taking $(\varepsilon_1, \varepsilon_2, \varepsilon_3) = (\varepsilon, \varepsilon, \varepsilon)$, we arrive in $R_1$ and we stop there. Therefore, $u_{1, 2, 3}(0)=\frac{1}{\lambda}\ell_1(0, a_1)$ cannot be lower than $V_{sc(1)}(0)$ that is a contradiction.

3.2) If $\sigma=23$ we take the triple $(\varepsilon_1, \varepsilon_2, \varepsilon_3) = (\varepsilon^2, \varepsilon, \varepsilon)$ and argue as in case 2).\\
We remark that if $\sigma=12$ or $13$ then considering $ (\varepsilon, \varepsilon, \varepsilon^2)$ and $ (\varepsilon, \varepsilon^2, \varepsilon)$ respectively we can conclude as in 3.1).



4) $f_1(0,a_1), f_2(0, a_2)=0, f_3(0, a_3)<0$. Also in this case we have different sub-cases according to the value of $\sigma$.

4.1) If $\sigma=123$ we take the triple $(\varepsilon, \varepsilon, \varepsilon)$ and conclude using Remark \ref{Remarkcontrolli} since $u_{1, 2, 3}(0)$ cannot be lower than a state constraints. Then as before a contradiction.

4.2) If $\sigma=23$ taking  the triple $(\varepsilon^2, \varepsilon, \varepsilon)$ we get $u_{2, 3}(0)= \frac{1}{\lambda}\ell_2(0, a_2)$ that is no lower than$ V_{sc(2)}(0)$, that is a contradiction.


Now we assume $\liminf V_{\varepsilon_1, \varepsilon_2, \varepsilon_3}(0) < V^{*}(0)$. Let $\delta > 0$ be such that, for arbitrarily small suitably chosen $(\varepsilon_1, \varepsilon_2, \varepsilon_3)$ , it is $V_{\varepsilon_1, \varepsilon_2, \varepsilon_3}(0) + \delta < V^{*}(0)$. A measurable control $\alpha$ which almost realizes the optimum (less then $\beta$) for $V_{\varepsilon_1, \varepsilon_2, \varepsilon_3}(0)$ must be such that there are infinitely many switching between all branches $R_{\varepsilon_1}, R_{\varepsilon_2}, R_{\varepsilon_3}$. Indeed, if it is not the case, then, for a least one branch, say $R_{\varepsilon_i}$, the trajectory definitely remains inside it. Hence, for small $(\varepsilon_1, \varepsilon_2, \varepsilon_3)$, $V_{\varepsilon_1, \varepsilon_2, \varepsilon_3}(0)$ is almost equal to $V_{sc(i)}(0)$, which is a contradiction. As in \cref{limitsogliefisse}, we can limit to consider a piecewise constant control that we call again $\alpha$. To prove that  $V_{\varepsilon_1, \varepsilon_2, \varepsilon_3}(0)$ cannot be less than $V^{*}(0) - \delta$ we proceed as in \cref{limitsogliefisse} considering $O(\max\{\varepsilon_1, \varepsilon_2, \varepsilon_3\})$ and $u_0(0)$ instead of $O(\varepsilon)$ and $u_{1, 2, 3}(0)$ respectively.\\
In conclusion we have $\liminf V_{\varepsilon_1, \varepsilon_2, \varepsilon_3}(0) = V^{*}(0)$.  The equations solved by $V_{\varepsilon_1, \varepsilon_2, \varepsilon_3}$ and by $V^{*}(0)$(\cref{RO,s7,s8} suitably modified) are the same in the interior of $R_i$ and the boundary datum converges to $V^{*}(0)$. Then, representing the solutions as the value functions of the corresponding optimal control problems, we get \cref{s11} and the uniform convergence.
\end{proof}
\begin{rem} 
As we show in the proof of \cref{limitesogliemobili}, we can restrict us to consider as thresholds $
(\varepsilon, \varepsilon, \varepsilon ), (\varepsilon, \varepsilon, \varepsilon^2 ), (\varepsilon, \varepsilon^2, \varepsilon ), (\varepsilon^2, \varepsilon, \varepsilon )$.
Hence, given the dynamics $f_1, f_2, f_3$ and the running costs $\ell_1, \ell_2, \ell_3$ satisfying the controllability assumptions exists a unique choice of previous triples such that
\begin{equation*}
V^*(x, i)=\liminf_{(\varepsilon_1, \varepsilon_2, \varepsilon_3)\rightarrow (0, 0, 0)} V_{\varepsilon_1, \varepsilon_2, \varepsilon_3}(x, i)= \lim_{(\cdot, \cdot, \cdot)\rightarrow (0, 0, 0)}V_{(\cdot, \cdot, \cdot)}(x, i) \ \forall\ (x, i) \in R_i.
\end{equation*}
We do not consider triples of the kind $(c_1\varepsilon, c_2\varepsilon, c_3\varepsilon)$, $c_1, c_2, c_3 \in \RR$ because they do not bring new possible optimal behaviours. Similarly, we do not consider triples as $(\varepsilon^2, \varepsilon^2, \varepsilon)$ because at the limit this would means to stay in $x = 0$ without
using the balance of the dynamics, which is physically meaningless.
\end{rem}
\begin{rem}
When the optimal strategy is to switch among all branches we
have that $\ell_i(0, a_i) = \ell_j(0, a_j)\ \forall i, j \in \{1, 2, 3\}, i\neq j$ and $V^*(x, i)=V(x, i)$, where $V$ is the value function \cref{funzionevaloresogliefisse} of the threefold junction problem with uniform thresholds.
\end{rem}
We introduce test functions $\psi: TR\to \RR$ such that $\psi \in C^1(TR)$ and on each branch $\psi_i: R_i\to \RR$ is such that  $\psi_i(x, i)=\psi_j(x, j)$ for every $i,j \in \{1, 2, 3\}, i\neq j$ when $x=0$.
\begin{definition}
\label{def:visco}
A continuous function $u: TR\to \RR$ is a viscosity subsolution of \cref{eq:HJBproblem3'} if for any $(x, i) \in TR$, for any test function $\psi$ as above such that $u-\psi$ has a local maximum at $(x, i)$, then
\begin{equation}
\begin{aligned}
& \lambda u(x, i) + H_i(x, \psi'_i(x, i))\leq 0 & & (x, i)\in int(R_i), \\
&  \min\left\{\lambda u(0, i)+H_i(0, \psi'_i(0, i)),\  i=1, 2, 3 \right\}\leq 0 & & x=0;
\end{aligned}
\end{equation}
A continuous function $u: TR\to \RR$ is a viscosity supersolution of \cref{eq:HJBproblem3'} if for any $(x, i) \in TR$, any $\psi \in C^1(TR)$ such that $u-\psi$ has a local minimum at $(x, i)$, then
\begin{equation}
\begin{aligned}
& \lambda u(x, i) + H_i(x, \psi'_i(x, i))\geq 0 & & x\in int(R_i), \\
&  \max\left\{\lambda u(0, i)+H_i(0, \psi'_i(0, i)),\ i=1, 2, 3 \right\}\geq 0  & & x=0.
\end{aligned}
\end{equation}
In particular, if $x=0$ then the local maximum/minimum may be considered
with respect to two of the three branches only.
\end{definition} 
Note the difference with \cref{defsubsup} where, for $x = 0$, the max/min is respect to all three branches.
\begin{theorem} Assume \cref{eq:Lip,eq:Controllability,eq:LLip}. The function $V^{*}$ is a viscosity solution and the maximal subsolution of the HJB problem, in the sense of \cref{def:visco}. 
\begin{equation}\label{eq:HJBproblem3'}
\begin{cases}
\lambda V + H_{1}(x, V') = 0 \ \ \text{in} \ int(R_1), \\
\lambda V + H_{2}(x, V') = 0 \ \ \text{in} \ int(R_2), \\
\lambda V + H_{3}(x, V') = 0 \ \ \text{in} \ int(R_3), \\
\min \left\{\lambda V + H_1, \lambda V + H_{2}, \lambda V + H_{3}\right\} \leq 0 \ \text{on} \ x=0, \\
\max \left\{\lambda V + H_1, \lambda V + H_{2}, \lambda V + H_{3}\right\} \geq 0 \ \text{on} \ x=0.
\end{cases}
\end{equation}
\end{theorem}
\begin{proof} 
By \cref{Approx3,limitesogliemobili}, $V^*$ satisfies the first three equations of \cref{eq:HJBproblem3'}.
For the other equations suppose, for example, that the optimal strategy is to switch only between the two branches $R_1$ and $R_2$. \\ Then $V^*= \lim_{(\varepsilon, \varepsilon, \varepsilon^2) \to (0, 0, 0)}V_{\varepsilon, \varepsilon, \varepsilon^2}=V_{1,2}$. If $V^*-\psi$ assumes its maximum or minimum in $x=0$ with respect to $R_1\cup R_2$, then, by the twofold junction problem, $V^*$ is a viscosity solution and the maximal subsolution of 
\begin{equation*}
\begin{cases}
\lambda V + H_{1}(x, V') = 0 \ \ \text{in} \ \ int(R_1),\\
\lambda V + H_{2}(x, V') = 0 \ \ \text{in} \ \ int(R_2), \\
\min \left\{\lambda V + H_1, \lambda V + H_{2}\right\} \leq 0 \ \text{on} \ x=0, \\
\max \left\{\lambda V + H_1, \lambda V + H_{2}\right\} \geq 0 \ \text{on} \ x=0.
\end{cases}
\end{equation*} 
If instead $V^*-\psi$ has a maximum point at $x=0$ with respect to $ R_1\cup R_3$ we prove that the $\min \left\{\lambda V + H_1, \lambda V + H_{2}\right\}$ is still lower or equal to zero. We consider two cases:

1) If the optimal behavior consists to reach $R_2$ and stay there, namely $V^*(x, 2)= V_{sc(2)}(x)$, and supposing that the cost to pay in $R_1$ to reach the junction is lower than the one in $R_3$, then $V^*=V_{1,2}$ on $R_1\cup R_2$. Now, since (by the assumption) $V^*-\psi$ has maximum point at $x=0$ locally with respect to the branch $R_3$, then $\psi_3(x, 3)\geq V^*(x, 3)$ for $x$ near to zero. The optimality of $V_{sc(2)}$ implies that
$V^*(\cdot, 3)\geq V_{sc(2)}(\cdot)=V^*(\cdot, 2)$ and hence $\psi_3(\cdot, 3)\geq V^*(\cdot, 2)$. Then gluing $\psi_3$ over $R_2$ we obtain that $V^*-\psi_3$ has a maximum point in $x=0$ locally with respect to $R_2$. Hence, $\min \left\{\lambda V + H_1, \lambda V + H_{2}\right\}\leq 0$.

2) If the optimal strategy is to switch between $R_1$ and $R_2$ and the maximum point at $x=0$ is still with respect to $R_1\cup R_3$ we conclude because $\psi_3(\cdot, 3)\geq V^*(\cdot, 2)$.

If $V^*=V_{1,2}$ and $V^*-\psi$ has a maximum point at $x=0$ with respect to $R_2\cup R_3$, with similar argument as before we conclude that $\min \left\{\lambda V + H_1, \lambda V + H_{2}\right\} \leq 0$.\\
In conclusion we have shown that the following condition hold:
exists a couple of indexes $(\bar{i}, \bar{j})$, fixed a priori, such that $V^*=V_{\bar{i}, \bar{j}}$ on $R_{\bar{i}}\cup R_{\bar{j}}$ and that for all $\psi \in C^1(TR)$ such that $V^*-\psi$ has the maximum point at $x=0$ with respect to any couple of edges, $\min\lbrace \lambda V + H_{\bar{i}}, \lambda V + H_{\bar{j}} \rbrace \leq 0$. From the latter condition follows that $\min\lbrace \lambda V + H_1, \lambda V + H_{2}, \lambda V + H_{3} \rbrace \leq 0$. Proceeding as before also for the fifth equation of \cref{eq:HJBproblem3'} we have that $V^*$ is a viscosity solution of \cref{eq:HJBproblem3'}.
Now, let $u$ be a continuous subsolution of \cref{eq:HJBproblem3'} satisfying the above condition with the same couple of indexes $(\bar{i}, \bar{j})$ that we suppose to be $(1, 2)$. Then 
\begin{equation}\label{cond su u}
V^*\geq u \ \ \text{on} \ \ R_1\cup R_2\ \Longrightarrow\ V^*(0)\geq u(0).
\end{equation}
Furthermore $V^*$ is a supersolution of the third equation of \cref{eq:HJBproblem3'}, $u$ is a subsolution of the same equation and hence by \cref{cond su u} follows $V^*\geq u$ on $R_3$. We can conclude that $V^*\geq u$ on $TR$ and hence it is the maximal subsolution of \cref{eq:HJBproblem3'}.
\end{proof}
Similarly $V^*= \lim_{(\varepsilon^2, \varepsilon, \varepsilon) \to (0, 0, 0)}V_{\varepsilon^2, \varepsilon, \varepsilon}=V_{2, 3}$, $V^*= \lim_{(\varepsilon, \varepsilon^2, \varepsilon) \to (0, 0, 0)}V_{\varepsilon, \varepsilon^2, \varepsilon}=V_{1, 3}$. 
\section*{Appendix}
	{\it Proof of the Proposition 4.1.}\\
		As in \cref{th:sistemapprox} in \cite{Ba} by virtue of the total controllability and Dynamic Programming (at least for small $\varepsilon$) we have 
		\begin{multline}\label{exittimevaluefunction}
			V_{\varepsilon_1, \varepsilon_2, \varepsilon_3}(x_0, i_0)= \inf_{\alpha \in {\cal A}}\Biggl\{ \int_{0}^{t_{(x_0, i_0)}(\alpha)}e^{-\lambda s}\ell_{i_0}(x(s), \alpha(s))ds \\ + e^{-\lambda t_{(x_0, i_0)}(\alpha)}V_{\varepsilon_1, \varepsilon_2, \varepsilon_3}(-\varepsilon_{i_{0^+}}, i_{0^+}) \Biggr\}.
		\end{multline}
		Namely, in each connected component $R_{\varepsilon_i}$, $V_{\varepsilon_1, \varepsilon_2, \varepsilon_3}$
		is the value function of the exit time problem from $R_{\varepsilon_i}$ with exit cost $V_{\varepsilon_1, \varepsilon_2, \varepsilon_3}(-\varepsilon_{i_{0^+}}, i_{0^+})$. Here $t_{(x_0, i_0)}(\alpha)$ is the first switching time, $i_{0^+}$ is the next value of the output $i_0$ and $\varepsilon_{i_{0^+}}$ the relative threshold. We get that the value function $V_{\varepsilon_1, \varepsilon_2, \varepsilon_3}$ is bounded and uniformly continuous on each of the three connected components of $ {\cal O}=R_{\varepsilon_1} \cup R_{\varepsilon_2} \cup R_{\varepsilon_3}$. Put together all these considerations, by a classical result on the boundary conditions in the viscosity sense follows that 
		$V_{\varepsilon_1, \varepsilon_2, \varepsilon_3}$ is a viscosity solution of the system \eqref{RO} on each branch with condition in viscosity sense. Regarding the uniqueness we prove that every solution of \eqref{RO} is a fixed point of a contraction map ${\cal G}: BC({\cal O})\to BC({\cal O})$, where $BC({\cal O})= BC(R_{\varepsilon_1})\times BC(R_{\varepsilon_2})\times BC(R_{\varepsilon_3})$ is the space of the bounded and continuous function on ${\cal O}$. Hence, by completeness arguments the uniqueness follows. We sketch how the contraction map ${\cal G}$ is constructed.
		For every $c\geq 0$ and for every $i_0 \in \lbrace 1, 2, 3\rbrace$, let $z_c^{(i_0)}$ be the solution of the correspondent Hamilton-Jacobi equation \eqref{RO}$_{i_0}$ with boundary datum $c$. Hence, for each $(\xi, \eta, \sigma)\in BC({\cal O})$ we define
		\begin{equation*}
			{\cal G}(\xi, \eta, \sigma):= \biggl(z_{\bigl( z_{\xi(-\varepsilon_2)}^{(2)}(\varepsilon_2) \bigr)}^{(1)}(\cdot), z_{\bigl( z_{\eta(-\varepsilon_3)}^{(3)}(\varepsilon_3) \bigr)}^{(2)}(\cdot), z_{\bigl( z_{\sigma(-\varepsilon_1)}^{(1)}(\varepsilon_1) \bigr)}^{(3)}(\cdot)\biggr).
		\end{equation*}
		This means that, for instance, the first component of ${\cal G}(\xi, \eta, \sigma)$ is the solution on the branch $R_{\varepsilon_1}$ with boundary datum equal to the value on $\varepsilon_2$ of the solution on the branch $R_{\varepsilon_2}$ with boundary datum equal to $\xi(-\varepsilon_2)$. Then for every $(\xi, \eta, \sigma), (\widehat{\xi}, \widehat{\eta}, \widehat{\sigma}) \in BC({\cal O})$, for the first component of ${\cal G}$ we have
		\begin{align*}
			& \Vert ({\cal G}(\xi, \eta, \sigma))_1 - ({\cal G}(\widehat{\xi}, \widehat{\eta}, \widehat{\sigma}))_1\Vert_{\infty} \leq \vert z_{\xi(-\varepsilon_2)}^{(2)}(\varepsilon_2) - z_{\widehat{\xi}(-\varepsilon_2)}^{(2)}(\varepsilon_2)\vert \\
			& \leq e^{\frac{-\lambda (2\varepsilon_2)}{M}}\vert \xi(-\varepsilon_2)- \widehat{\xi}(-\varepsilon_2)\vert \leq e^{\frac{-\lambda (2\varepsilon_2)}{M}} \Vert \xi - \widehat{\xi}\Vert_{\infty},
		\end{align*}
		with $M$ the bound of the dynamics $f_i$. A similar inequality holds for the others components of ${\cal G}$. Since $\lambda>0$ we get the conclusion.
%
%
%
\bibliographystyle{siamplain}
\bibliography{references}
\end{document}


\maketitle

\section{A detailed example}

Here we include some equations and theorem-like environments to show
how these are labeled in a supplement and can be referenced from the
main text.
Consider the following equation:
\begin{equation}
  \label{eq:suppa}
  a^2 + b^2 = c^2.
\end{equation}
You can also reference equations such as \cref{eq:matrices,eq:bb} 
from the main article in this supplement.

\lipsum[100-101]

\begin{theorem}
  An example theorem.
\end{theorem}

\lipsum[102]
 
\begin{lemma}
  An example lemma.
\end{lemma}

\lipsum[103-105]

Here is an example citation: \cite{KoMa14}.

\section[Proof of Thm]{Proof of \cref{thm:bigthm}}
\label{sec:proof}

\lipsum[106-114]

\section{Additional experimental results}
\Cref{tab:foo} shows additional
supporting evidence. 

\begin{table}[htbp]
  \caption{Example table}
  \label{tab:foo}
  \centering
  \begin{tabular}{|c|c|c|} \hline
   Species & \bf Mean & \bf Std.~Dev. \\ \hline
    1 & 3.4 & 1.2 \\
    2 & 5.4 & 0.6 \\ \hline
  \end{tabular}
\end{table}

\bibliographystyle{siamplain.bst}
\bibliography{references}